\documentclass{amsart}

\usepackage{amsmath,amscd,amssymb,amsthm}
\usepackage{enumerate,mathrsfs}
\usepackage{geometry}
\usepackage{hyperref}
\usepackage{url}
\usepackage{ifthen}
\usepackage[all]{xy} \xyoption{2cell} \UseAllTwocells

\newcommand{\bibfilename}{/home/micromath/results/b20140922unibib/unibib}
\makeatletter
\let\c@equation\c@subsection

\makeatother

\swapnumbers
\theoremstyle{plain}
\newtheorem{thm_}[subsection]{Theorem}
\newtheorem{lemma_}[subsection]{Lemma}
\newtheorem{prop_}[subsection]{Proposition}
\newtheorem{cor_}[subsection]{Corollary}
\newtheorem{eg_}[subsection]{Example}
\newtheorem{con_}[subsection]{Conjecture}
\newtheorem*{cons_}{Conjecture}

\theoremstyle{definition}
\newtheorem{thmu_}[subsection]{Theorem}
\newtheorem*{thmus_}{Theorem}
\newtheorem{propu_}[subsection]{Proposition}
\newtheorem*{propus_}{Proposition}
\newtheorem{coru_}[subsection]{Corollary}
\newtheorem{lemu_}[subsection]{Lemma}
\newtheorem*{lemus_}{Lemma}
\newtheorem{egu_}[subsection]{Example}
\newtheorem*{egus_}{Example}
\newtheorem{def_}[subsection]{Definition}
\newtheorem*{defs_}{Definition}
\newtheorem{rk_}[subsection]{Remark}

\newcommand{\thm}[1]{\begin{thm_}#1\end{thm_}}
\newcommand{\thmu}[1]{\begin{thmu_}#1\end{thmu_}}

\newcommand{\lemm}[1]{\begin{lemma_}#1\end{lemma_}}

\newcommand{\egu}[1]{\begin{egu_}#1\end{egu_}}

\newcommand{\prop}[1]{\begin{prop_}#1\end{prop_}}

\newcommand{\defi}[1]{\begin{def_}#1\end{def_}}

\newcommand{\rk}[1]{\begin{rk_}#1\end{rk_}}
\newcommand{\cor}[1]{\begin{cor_}#1\end{cor_}}
\newcommand{\coru}[1]{\begin{coru_}#1\end{coru_}}

\newcommand{\pf}[1]{\begin{proof}#1\end{proof}}

\DeclareMathOperator{\Hom}{Hom}

\DeclareMathOperator{\Spec}{Spec}

\DeclareMathOperator{\Aut}{Aut}
\DeclareMathOperator{\im}{im}
\DeclareMathOperator{\Ob}{Ob}
\DeclareMathOperator{\Br}{Br}

\DeclareMathOperator{\id}{id}
\DeclareMathOperator{\ad}{ad}
\newcommand{\bA}{\mathbb A}%
\newcommand{\CC}{\mathbb C}
\newcommand{\FF}{\mathbb F}

\newcommand{\II}{\mathbb I}

\newcommand{\ZZ}{\mathbb Z}
\newcommand{\bfA}{\mathbf A}%
\newcommand{\bfG}{\mathbf G}%

\newcommand{\sE}{\mathscr E}
\newcommand{\sF}{\mathscr F}
\newcommand{\sG}{\mathscr G}
\newcommand{\sH}{\mathscr H}

\newcommand{\sU}{\mathscr U}
\newcommand{\sX}{\mathscr X}
\newcommand{\sY}{\mathscr Y}
\newcommand{\sZ}{\mathscr Z}
\newcommand{\cA}{\mathcal A}%
\newcommand{\cC}{\mathcal C}%
\newcommand{\cD}{\mathcal D}
\newcommand{\cG}{\mathcal G}

\newcommand{\cO}{\mathcal O}%

\newcommand{\cR}{\mathcal R}
\newcommand{\cS}{\mathcal S}
\newcommand{\cT}{\mathcal T}
\newcommand{\cU}{\mathcal U}
\newcommand{\cV}{\mathcal V}
\newcommand{\cX}{\mathcal X}
\newcommand{\cY}{\mathcal Y}
\newcommand{\cZ}{\mathcal Z}%
\newcommand{\sfH}{\mathsf H}
\newcommand{\s}{\sigma}
\newcommand{\La}{\Lambda}
\newcommand{\tm}{\times}%

\newcommand{\ul}{\underline}
\newcommand{\ra}{\rightarrow}

\newcommand{\Ra}{\Rightarrow}
\newcommand{\lra}{\longrightarrow}

\newcommand{\hra}{\hookrightarrow}

\newcommand{\mpt}{\mapsto}

\newcommand{\os}[2]{\overset{#1}{#2}}
\newcommand{\tu}[1]{\text{\upshape #1}}
\newcommand{\PPt}{{Pt}}
\newcommand{\MMor}{{Mor}}
\newcommand{\SSet}{{Set}}

\newcommand{\SSch}{{Sch}}
\newcommand{\AAb}{{Ab}}
\newcommand{\MMod}{{Mod}}
\newcommand{\CCat}{{\mathbb C\tu{at}}}
\newcommand{\TTps}{{\mathbb T\tu{ps}}}

\newcommand{\FFib}{{\mathbb F\tu{ib}}}
\newcommand{\FFibS}{{Fib}_\SSet}
\newcommand{\FFibP}{{\FFib_{PSh}}}
\newcommand{\FFibR}{{\FFib_{Rep}}}
\newcommand{\SStk}{{\mathbb S\tu{tk}}}
\newcommand{\SStkS}{{Stk}_\SSet}
\newcommand{\SStkP}{{\SStk_{Sh}}}
\newcommand{\PPSh}{{PSh}}
\newcommand{\SSh}{{Sh}}
\newcommand{\SSGrp}{{SGrp}}

\newcommand{\et}{\tu{\'et}}

\newcommand{\fppf}{\tu{fppf}}

\newcommand{\Tors}{\mathsf{Tors}}
\newcommand{\Gerb}{\mathsf{Gerb}}
\newcommand{\eq}[1]{\begin{equation}#1\end{equation}}
\newcommand{\eqn}[1]{\begin{equation*}#1\end{equation*}}
\newcommand{\ga}[1]{\begin{gather}#1\end{gather}}
\newcommand{\gan}[1]{\begin{gather*}#1\end{gather*}}
\newcommand{\al}[1]{\begin{align}#1\end{align}}
\newcommand{\aln}[1]{\begin{align*}#1\end{align*}}

\newcommand{\enmt}[1]{\begin{enumerate}#1\end{enumerate}}

\newcommand{\sptag}[1]{\href{https://stacks.math.columbia.edu/tag/#1}{#1}}
\newcommand{\citespt}[3]{%
\ifthenelse{\equal{#2}{}}{
\cite[Tag \sptag{#1}]{stacks-project}}{%
\ifthenelse{\equal{#3}{}}{
\cite[Tags \sptag{#1} and \sptag{#2}]{stacks-project}}{
\cite[Tags \sptag{#1}, \sptag{#2} and \sptag{#3}]{stacks-project}}}}
\newcommand{\citesp}[1]{\citespt{#1}{}{}}
\newcommand{\citespp}[2]{\citespt{#1}{#2}{}}
\newcommand{\citesppp}[3]{\citespt{#1}{#2}{#3}}

\newcommand{\ob}{\tu{ob}}
\newcommand{\deltaob}{{\delta, \ob}}
\newcommand{\cdesc}{\tu{cdesc}}
\newcommand{\desc}{\tu{desc}}
\newcommand{\descob}{{\desc, \ob}}
\newcommand{\dcdesc}{{\cdesc, \cdesc}}
\newcommand{\etBr}{{\et, {\Br}}}
\newcommand{\sdesc}{{2\tu{-}\desc}}
\newcommand{\sdescob}{{\sdesc, \ob}}

\newcommand{\hdesc}[2]{{#1\tu{-}\desc{\ifthenelse{\equal{#2}{}}{}{_#2}}}}
\newcommand{\sndesc}{{2\tu{-ndesc}}}

\newcommand{\SSht}{{\SSh_\tau}}
\newcommand{\TTpsSt}{{\TTps/\SSh(\cS_\tau)}}
\newcommand{\St}[1]{\SSh(#1_\tau)}
\newcommand{\XA}{X(\bfA_k)}
\newcommand{\YA}{Y(\bfA_k)}
\newcommand{\Xk}{X(k)}

\newcommand{\cXT}{\cX(\cT)}
\newcommand{\cXA}{\cX(\cA)}
\newcommand{\cXS}{\cX(\cS)}
\newcommand{\cYA}{\cY(\cA)}

\newcommand{\cYT}{\cY(\cT)}
\newcommand{\cUA}{\cU(\cA)}
\newcommand{\cVA}{\cV(\cA)}
\newcommand{\DM}{Deligne-Mumford }

\begin{document}
\title[Second descent and gerbes]{The second descent obstruction and gerbes}
\author[C. Lv]{Chang Lv}
\address{State Key Laboratory of Information Security\\
Institute of Information Engineering\\
Chinese Academy of Sciences\\
Beijing 100093, P.R. China}
\email{lvchang@amss.ac.cn}
\date{\today}
\thanks{This work was supported by
 National Natural Science Foundation of China (Grant No. 11701552).
}
\begin{abstract}
We extend the notion of rational points and cohomological obstructions
 on varieties to  categories fibred in groupoids.
We also establish the generalized theory of descent by torsors.
Then we  interpret the obstruction given by the second cohomology of abelian
 sheaves in terms of categorical points of gerbes, analogue to descent by
 torsors.
As an application, we construct some composite obstructions not larger than
 the descent obstruction.
We also propose some new kinds of obstructions including the  derived
 obstruction, which has
 good behavior under a product and is also not larger than the descent
 obstruction.
\end{abstract}
\maketitle

\setcounter{tocdepth}{1}
\tableofcontents

\section{Introduction}\label{sec_intro}
The existence of rational points is a fundamental arithmetic property of
 varieties over number fields.
One of the most common methods to consider is the local-global principle.
When it fails, one considers various obstructions to it, such as the
 Brauer-Manin obstruction and the descent obstruction.
In this section, we review the classical descent set and descent by torsors.
Then we introduce our results on the generalizations of points,
 obstructions, descent by gerbes, and derived obstructions.

\subsection{Rational points on varieties and the local-global principle}
Let $X$  be a variety over a number field $k$, we write $X(-)$ for the functor
 $\Hom_k(-, X)$. The \emph{rational points} of $X$ is the set $\Xk$, which is
 contained in the \emph{ad\`elic points} $\XA$.
We say that the \emph{local-global principle} holds if $\XA\neq\emptyset$
 implies $\Xk\neq\emptyset$.

\subsection{The $F$-obstruction to the local-global principle} \label{Fob_c}
Let $F: (\SSch/k)^\circ\ra \SSet$ be a contravariant
 functor from the category of $k$-schemes
 to the category of sets.
For a $k$-scheme $T$ and $A\in F(X)$, the \emph{evaluation} of
 a \emph{$T$-point} $T\os{x}{\ra}X \in X(T)$ at
 $A$ is defined to be  the image  of  $A$ under the pull-back map
 $F(x): F(X)\ra F(T)$  induced by $x$, denoted by $A(x)$.
We have an obvious commutative diagram
\eqn{\xymatrix{
\Xk\ar[r]\ar[d]^-{A(-)} &\XA\ar[d]^-{A(-)} \\
F(k)\ar[r] &F(\bfA_k)
}}
Define the \emph{obstruction $($set$)$ given by $A$} to be
  $\XA^A=\{x\in \XA\mid A(x)\in \im(F(k)\ra F(\bfA_k)\}$. Then we
 have $\emptyset\subseteq \Xk\subseteq \XA^A\subseteq \XA$, giving
 a constraint on the locus of  rational points in ad\`elic points.
Putting
\eqn{
\XA^F=\XA^{F(X)}=\bigcap_{A\in F(X)} \XA^A,
}
 called the \emph{$F$-set}, which also yields the inclusion
 $\Xk\subseteq \XA^F\subseteq\XA^A\subseteq\XA$.
See Poonen \cite[8.1.1]{poonen17rational}.

\defi{
We say that the $F$-obstruction to the local-global principle
 is the \emph{only one} if
 $\XA^F\neq\emptyset$ implies $\Xk\neq\emptyset$.

We say that \emph{there is} an
 $F$-obstruction to the local-global principle if $\XA\neq\emptyset$ but
 $\XA^F=\emptyset$ (a \emph{priori} $\Xk=\emptyset$).
}

\subsection{Examples} \label{eg_cohob}
We are mainly interested in \emph{cohomological obstructions}, namely, ones
 where $F$ is taken to be a cohomological functor.

Let $F=\Br=H^2_{\text{\'et}}(-, \bfG_m)$ be
 the cohomological Brauer-Grothendieck group (c.f. Grothendieck
 \cite{grothendieck95braueri, grothendieck95brauerii, grothendieck95braueriii}).
Thus we obtain $X(\bA_k)^{\Br}$, the Brauer-Manin set (see, e.g. Skorobogatov
 \cite{torsor}).

Another example is that $F=\check H_\fppf^1(-, G)$, the obstruction given by
 the first \v Cech cohomology, where $G$ is an affine $k$-group. In general,
 $\check H_\fppf^1(X, G)$ is a pointed set, which is isomorphic to
 $H_\fppf^1(X, G)$ if $G$ is commutative, and further  to $H_\et^1(X, G)$ since
 $G$ is smooth over the number field $k$.
The classical \emph{descent set} is given by
\eq{ \label{eq_cdesc}
\XA^\cdesc=
 \bigcap_{\text{all affine $k$-group $G$}} \XA^{\check H_\fppf^1(X, G)}.
}
\rk{
In Section \ref{desc} we shall generalize the classical descent set  from
 varieties to fibred categories.
To avoid confusion, here we use $\cdesc$ instead of $\desc$ to indicate that
 it is the ``classical" descent set.
}
The descent theory was established by Colliot-Th{\'e}l{\`e}ne  and Sansuc
 \cite{colliot1987descente} for tori  and Skorobogatov
 \cite{skorobogatov99beyond} for groups of multiplicative type.
Harari \cite{harari02groupes}, Harari
 and Skorobogatov \cite{hs02non-abelian,
 hs05non-abelian} studied the descent obstruction for general algebraic groups
 and compared it with the Brauer-Manin obstruction.
One of the results is the well-known inclusion
 $\XA^\cdesc\subseteq \XA^{\Br}$ for regular, quasi-projective $k$-variety $X$
 (see, e.g., \cite[Prop. 8.5.3]{poonen17rational}).

The classical descent by torsors says that
\eqn{
\XA^f=\bigcup_{\s\in \check H^1(k, G)} f^\s(Y^\s(\bfA_k))
}
 where $f: Y\os{G}{\ra} X$ is a $G$-torsor over $X$,
 representing the class $[Y]\in \check H_\fppf^1(X, G)$,
 and $f^\s: Y^\s\os{G^\s}{\ra} X$ is the twisted torsor.
Then we may define various of composite (or iterated) obstruction sets
 $\XA^{\cdesc, \ob}$ between $\Xk$ and $\XA^\cdesc$, such as
\eqn{
\XA^\dcdesc=
 \bigcap_{\text{all affine $G$ and all torsor } f: Y\os{G}{\ra} X}
 \bigcup_{\s\in \check H^1(k, G)} f^\s(Y^\s(\bfA_k)^\cdesc)
}
 and
\eqn{
\XA^\etBr=
 \bigcap_{\text{all finite $G$ and all torsor } f: Y\os{G}{\ra} X}
 \bigcup_{\s\in \check H^1(k, G)} f^\s(Y^\s(\bfA_k)^{\Br}).
}
By works of Stoll, Skorobogatov, Demarche, Poonen, Xu and Cao,
 it is known that for $X$ being a
 smooth, quasi-projective, geometrically integral
 $k$-variety,  $\XA^\cdesc=\XA^\etBr$ \cite{stoll07finite,
 skorobogatov09descent, demarche09obstruction, poonen10insufficiency,
 cdx19comparing} and $\XA^\dcdesc=\XA^\cdesc$ \cite{cao20sous} and no smaller
 obstruction set than the descent set is discovered.

Next we begin to introduce main results of this text.

\subsection{The second descent and composite obstructions}
For any commutative $k$-group $G$ we have \cite{harari02groupes}
 $\XA^{\Br}\subseteq \XA^{H_\et^2(X, G)}$ for smooth geometrically integral
 $k$-variety $X$.
But we want to obtain subsets not larger than the descent set.
Define the \emph{second descent set} (see \ref{defi_sdesc})
\eqn{
\XA^\sdesc=\bigcap_{\text{all abelian sheaf $\sG$ on $(\SSch/k)_\et$}}
 \XA^{H_\et^2(X, \sG)}.
}
Since  each commutative $k$-group  $G$ represents an abelian sheaf on the big
 \'etale site $(\SSch/k)_\et$,
 it is clear that $\XA^\sdesc\subseteq \XA^{H_\et^2(X, G)}$.
In particular, $\XA^\sdesc\subseteq \XA^{\Br}$.
However we do not know whether $\XA^\sdesc\subseteq \XA^\cdesc$ holds.

The key point to go further is to use points of
 gerbes to interpret the obstruction given
 by second cohomology, which is analogue to classical descent by torsors.
Recall that a gerbe is a kind of categories fibred in groupoids.
In fact, most of these obstructions on $\XA$ mentioned above
 can be generalized to that on $\cXA$,
 where $\cX$ and $\cA$ are categories over $\cS$ fibred in groupoids, $\cS
 =\SSch/S$ is the category of $S$-schemes.
To go back to the classical case, one sets
 $S=\Spec k$, $\cX=X$ a variety over $k$  and take
 $\cA\ra \cS$ to be $\Spec\bfA_k\ra \Spec k$ induced by the
 inclusion $k\subset \bfA_k$, and then we have $\cXA=\XA$ and
 $\cXS=\Xk$.
See Section \ref{pts+obs}.

Assume further that $\cA$ (resp. $\cX$) is a stack over $\cS_\et$ fibred in
 sets (resp. groupoids).
Then we show that (descent by gerbes)
\eqn{
\cXA^f=\bigcup_{\s\in H_\et^2(\cS, \sG)} f^\s(\cY^\s(\cA))
}
 where $f: \cY\os{\sG}{\ra} \cX_\et$ is any fixed gerbe over $\cX_\et$ bounded
 by an abelian sheaf $\sG$ on $\cS_\et$,
 representing the class $[\cY]\in H_\et^2(\cX, \sG)$. See \ref{thm_dbg}.
Then we have
\eqn{
\cXA^\sdesc=
 \bigcap_{\text{all $\sG$ on $\cS_\et$ and all gerbes }
 f: \cY\os{\sG}{\ra} \cX_\et}
 \bigcup_{\s\in H_\et^2(\cS, \sG)} f^\s(\cY^\s(\cA)).
}
This leads to the composite obstruction
\eqn{
\cXA^\sdescob=
 \bigcap_{\text{all $\sG$ on $\cS_\et$ and all gerbes }
 f: \cY\os{\sG}{\ra} \cX_\et}
 \bigcup_{\s\in H_\et^2(\cS, \sG)} f^\s(\cY^\s(\cA)^\ob).
}
In general, taking into account the generalized descent by torsors developed in
 Section \ref{desc},
 we shall  show in \ref{prop_cXAH1} that $\XA^\desc\subseteq \XA^\cdesc$,
 and in \ref{cor_comp_ob} that,
\eq{ \label{eq_c1}
\cXS\subseteq \dots\subseteq \cXA^{\delta,\delta,\delta, \ob}\subseteq
 \cXA^{\delta,\delta, \ob}\subseteq \cXA^{\delta,\ob}\subseteq
 \cXA^\delta\cap \cXA^\ob\subseteq \cXA,
}
 where $\delta\in\{\desc,\sdesc\}$, and if $\delta=\desc$, we further assume
 $\cX$ is a stack in setoid, and where
 $\ob$ is an \emph{functorial} obstruction map (see \ref{defi_func_map}),
 such as $\desc$, $\sdesc$ and $\SSht$ below.

\subsection{The derived obstructions}
Finally  we define \emph{derived obstruction} $\cXA^\SSht$ by $1$-morphisms of
 topoi (see \ref{defi_do}), where $\tau\in\{\et, \fppf\}$.
With the corresponding $\tau$ chosen,
 this  obstruction is not larger than both descent obstruction and
 second descent obstruction, i.e.,
\ga{
\XA^{\SSh_\fppf}\subseteq \XA^\desc\subseteq\XA^\cdesc, \label{eq_c2} \\
\XA^{\SSh_\et}\subseteq \XA^\sdesc. \nonumber
}
See \ref{prop_do_cmp}.
We also show in \ref{thm_do_prod} that,
 under mild assumptions, it has good behavior under a product.
In particular, in the classical case, we have
\eqn{
(X\tm_k Y)(\bfA_k)^\SSht=\XA^\SSht\tm \YA^\SSht.
}

\subsection{Obstructions not larger than the classical descent obstruction}
From \eqref{eq_c1} and \eqref{eq_c2}, we obtain some obstructions contained in
 $\XA^\cdesc$:
\gan{
\XA^{\sdesc,\dots,\sdesc,\desc}\subseteq
 \XA^{\sdesc,\desc}\subseteq\XA^\desc\subseteq\XA^\cdesc, \\
\XA^{\delta,\dots,\delta,\ob}\subseteq
 \XA^{\delta,\ob}\subseteq\XA^\ob
 \subseteq\XA^\desc\subseteq\XA^\cdesc,
}
where $\ob\in\{\desc, \SSh_\fppf\}$ and $\delta$ as in \eqref{eq_c1}.

\subsection{} The paper is organized as follows.
Section \ref{pts+obs} gives a framework
 that extends the classical notion of rational points and
 cohomological obstructions on varieties to categories fibred in groupoids.
In Section \ref{desc} we reformulate the classical descent theory in the
 context of categories fibred in groupoids.
Section \ref{sdesc} is devoted to define the abelian second  descent
 obstructions on categories fibred in groupoids and formulate the idea of
 descent by gerbes.
We also define the non-abelian second descent obstruction.
Then in Section \ref{compob} we use results in previous two sections to
 give some composite obstructions not larger than the descent obstruction.
Finally, we propose  the higher descent obstruction
 and derived obstruction in Section \ref{hdescder}, and shows that
 the later has
 good behavior under a product and is also not larger than the descent
 obstruction.

\section{Generalizations on classical points and obstructions} \label{pts+obs}
In this section we extend the classical notion of rational points and
 cohomological obstructions on varieties to categories fibred in groupoids.

\subsection{Fibred categories} \label{fib}
If there is  no other confusion, we write $\cC\subseteq\cD$ (resp. $V\in\cD$)
 if $\cC$ is a full subcategory of $\cD$ (resp. $V\in\Ob(\cD)$).
If no extra explanation, \emph{fibred categories} mean categories fibred in
 groupoids.
Let $\cC$ be a category.
Then we have $\FFib/\cC$, the $2$-category of  fibred categories
 over $\cC$ \citesp{02XS}.
Recall that in $\FFib/\cC$,
 all $2$-morphisms in $\FFib/\cC$ are $2$-isomorphisms and
 $2$-fibred products exist.
Let $\cT'$ and $\cT$ be two fibred categories over $\cC$.
Then the category $\MMor_{\FFib/\cC}(\cT', \cT)$ in $\FFib/\cC$
 consists of, by definition, $1$-morphisms $\cT'\ra\cT$ as objects and
 $2$-morphisms between them as morphisms.
Thus it is a groupoid.
We also denote by $\CCat$ (resp. $\CCat/\cC$) the $(2, 1)$-category of
 categories (resp. categories over $\cC$)
 \citespp{003I}{003Y}.
We recall some basic definitions and properties.
\defi{ \label{defi_qs_equiv}
Let $\xymatrix{\cY\ar@<1mm>[r]^f &\cX\ar@<1mm>[l]^g}$ be two $1$-morphisms in
 a $2$-category $\CC$.
We say that $g$ is a \emph{quasi-section} (resp. \emph{section}) of $f$ if the
 diagram
\eqn{\xymatrix{
 &\cX\ar[dl]_-g \ar[d]^-\id \\
\cY\ar[r]^-f &\cX
}}
 $2$-commutes (resp. commutes) in $\CC$, that is,
 $f\circ g$ is $2$-isomorphic (resp. equal) to the identity in $\CC$.

We say that $f$ and $g$ define an \emph{equivalence} (resp. \emph{isomorphism})
 if they are quasi-sections (resp. sections) of each other.
If such $f$ and $g$ exist, $\cY$ and $\cX$ are \emph{equivalent}
 (resp. \emph{isomorphic}), denoted by $\cY\os{\approx}{\ra} \cX$
 (resp.  $\cY\os{\sim}{\ra} \cX$).
Note that if $\CC=\FFib/\cC$, an $1$-morphism $\cY\ra\cX$  in
 $\CCat/\cC$ that is an equivalence in $\CCat$
 is automatically an equivalence in $\FFib/\cC$ \citesp{003Z}.
}

\lemm{ \label{lemm_emb}
For an arbitrary category $\cC$, we have the following functors
\eq{ \label{eq_emb}
\xymatrix{
\cC \ar[r]^-{\epsilon_\cC} &\PPSh(\cC)\ar[r]^-{\eta_\cC}_-{\approx}
 &\FFibS/\cC\ \ar@{^(->}[r] &\FFib/\cC
}}
 where $\epsilon_\cC(U)=h_U=\MMor_\cC(-, U)$ is the Yoneda embedding
 $($presheaves in whose essential image are called \emph{representable}$)$,
 $\PPSh(\cC)$ is the category of presheaves on $\cC$,
 $\eta_\cC$ is an equivalence,
 $\FFibS/\cC$ is the category of categories over $\cC$ fibred in sets
 $($identified as a full sub $2$-category of $\FFib/\cC$ by the hook arrow$)$,
 and that:
\enmt{[\upshape (i)]
\item \label{it_eta}
 For $\sF\in\PPSh(\cC)$, the category $\eta_\cC(\sF)=\cC_\sF$ is defined as
  follows.
 The objects are pairs $(U, t)$ where $U\in\cC$ and $t\in\sF(U)$, and
 \eqn{
 \MMor_{\cC_\sF}((U',t'), (U,t))=\{f\in\MMor_\cC(U',U)\mid \sF(f)(t)=t'\},
 }
  and the structure morphism $\cC_\sF\ra \cC$ is the functor of
  projection to the first  component.
 The quasi-inverse of $\eta_\cC$ is defined to be, for $\cT\in\FFibS/\cC$,
  $\eta_\cC^{-1}(\cT)(U)
  =\cT_U$, the fiber category over $U$,
  and $\eta_\cC^{-1}(\cT)(f)=f^*$ where $f: U'\ra U$ and
  $f^*: \cT_U\ra\cT_{U'}$ is the
  pull-back functor in $\cT$.
\item \label{it_fib_psh}
 Let $\cT\in\FFib/\cC$.
 Then $\cT$ is  equivalent to $\eta_\cC(\sF)$ for some $\sF$
  $($unique up to an isomorphism$)$
  if and only if
  $\cT$ is fibred in  setoids.
 If this is the case,  we call $\cT$ \emph{is represented by presheaf},
  or \emph{is represented by $\sF$}.
 Denoted by $\FFibP/\cC$ the full sub $2$-category of $\FFib/\cC$ consisting
  categories over $\cC$ fibred in setoids.
\item \label{it_fib_rep}
 Let $\cT\in\FFib/\cC$.
 Then $\cT$ is  equivalent to $\eta_\cC\epsilon_\cC(U)$ for some $U$
  $($unique up to an isomorphism$)$
  if and only if $\cT$ is
  fibred in setoids and that the presheaf $U\mpt \Ob(\cT_U)_{/\cong}$ is
  representable.
 If this is the case, we call $\cT$  \emph{representable},
  or \emph{is represented by $U$}.
 Denoted by $\FFibR/\cC$ the full sub $2$-category of $\FFib/\cC$ consisting
  representable objects.
}}
\pf{
See \citesppp{02Y2}{0045}{02Y3}.
}
\rk{ \label{rk_emb}
Here are some remarks related to $\FFibS/\cC$.
\enmt{[\upshape (i)]
\item \label{it_U_Fib/S}
 If $U\in\cC$, then in \eqref{eq_emb}  the corresponding presheaf of sets
  $h_U$ represents an object $\eta_\cC\epsilon_\cC(U)$ in
  $\FFibS/\cC\hra\FFib/\cC$,
  it is the localized category $\cC/U$ \citesp{0044},
  still denoted by $U$.
\item \label{it_qs_s}
 In \ref{defi_qs_equiv},  let $\CC=\FFib/\cC$,
  if in addition $\cX\in\FFibS/\cC$, i.e., there is an
  isomorphism   $\cX\os{\sim}{\ra} \eta_\cC(\sF)$ for some $\sF\in\PPSh(\cC)$,
  then all quasi-sections of $f$ are sections of $f$ since
  $\MMor_{\FFib/\cC}(\cX, \cX)$ is a discrete category.
}}

\subsection{Pull-back of a fibred category} \label{pb_fib}
Let $u: \cT'\ra\cT$ be a $1$-morphism in $\FFib/\cC$ and
 $f: \cY\ra \cT\in\FFib/\cT$.
In particular, $\cY\in\FFib/\cC$ via compositing
 with  the structure morphism $\cT\ra\cC$  \citesp{09WW}.
\defi{ \label{defi_pb_fib}
Notations as in \ref{pb_fib},
 consider the base change of $\cY$ by $u$, i.e., the $2$-fibred product
 $u^*\cY=\cY\tm_\cT \cT'$
 indicated in the following $2$-cartesian diagram in $\FFib/\cC$
\eqn{\xymatrix{
u^*\cY=\cY\tm_\cT \cT'\ar[r]^-{f'}\ar[d]^-{u'} &\cT'\ar[d]^-u \\
\cY\ar[r]^-f &\cT
}}
Then $f': u^*\cY\ra \cT'$ is unique up to an equivalence and there is a choice
 making it a fibred category over $\cT'$ \citespp{003Q}{06N7}.
We shall always assume that $f': u^*\cY\ra\cT'$ is a fibred category, and
 call it the \emph{pull-back} of $\cY$ under $u$.
}
\lemm{ \label{lemm_pb_fib}
With notations in {\upshape \ref{defi_pb_fib}},
 assume further that $\cT', \cT\in\FFibP/\cC$.
Let $\cY'=u^*\cY$,  $t'\in\cT'$ and $t=u(t')\in\cT$, then $u$ induces an
 equivalence $\cY'_{t'}\os{\approx}{\ra}\cY_t$,
 where $\cY'_{t'}$ and $\cY_t$ are fiber categories.
}
\pf{
By \ref{lemm_emb} \eqref{it_fib_psh}, we may assume that
 $\cT', \cT\in\FFibS/\cS$.
Consider the following $2$-commutative diagram in $\FFib/\cC$
\eqn{\xymatrix{
\PPt\tm_{\cT'} \cY' \ar[r]\ar[d] &\cY'\ar[r]^-{u'}\ar[d]^-{f'} &\cY\ar[d]^-f \\
\PPt\ar[r]^-{t'}\ar@/_1pc/[rr]_-t &\cT'\ar[r]^-u &\cT
}}
 where $\PPt=\{*\}$ is a singleton category, and the arrow labeled with
 $t'$ (resp. $t$) is the functor sending $*$ to $t'$ (resp. $t$).
We see that all squares are $2$-cartesian and hence there is an equivalence
 $\PPt\tm_{\cT'} \cY'\os{\approx}{\ra} \PPt\tm_\cT \cY$.

Note that for any $\cR\in\FFib/\cC$, $\MMor_{\FFib/\cC}(\cR, \cT')$
 and $\MMor_{\FFib/\cC}(\cR, \cT)$ are
 discrete categories since $\cT', \cT\in\FFibS/\cS$.
Thus by the universal property of $2$-fibred product,  we have equivalences
 $\cY_t\os{\approx}{\ra} \PPt\tm_\cT \cY$ and
 $\cY'_{t'}\os{\approx}{\ra} \PPt\tm_{\cT'} \cY'$.
It follows that $\cY'_{t'}\os{\approx}{\ra} \cY_t$ and one checks that this is
 induced by $u'$, up to equivalences.
}

\defi{[Points]
For any $\cX, \cT\in\FFib/\cC$, we call  $\cXT=\MMor_{\FFib/\cC}(\cT, \cX)$
 the category of \emph{$\cT$-points in $\cX$}.
In particular, $\cX(\cC)$ is called \emph{rational} points of $\cX$.
}

\subsection{} \label{X(A)_X(S)} From now on, we work with $\FFib/\cS$, where
 $\cS=\SSch/S$ is the category of $S$-schemes, and $S$ is a fixed scheme.
We also fix a  fibred category $q: \cA\ra \cS$ in $\FFib/\cS$.
Let $p: \cX\ra \cS$ be another fibred category.
By abusing notations, we also use $\cXS$ for the essential image of the natural
 functor $\MMor_{\FFib/\cS}(q, \cX): \cXS\ra \cXA$ induced by $q$.
For example, we have $\emptyset\subseteq \cXS\subseteq \cXA$ (for the meaning
 of inclusions, see \ref{fib}).
Then we are interested  in full subcategories between $\cXS\subseteq \cXA$ and
 the ``local-global principle", which means that  $\cXA\neq\emptyset$ implies
 $\cXS\neq\emptyset$.
\egu{
Trivially we have $\cS(\cS)=\cS(\cA)=\{q\}$.
}

\rk{
To recover the classical problem of the local-global principle to rational
 points on varieties over a field $k$, let $\cX$ (resp. $\cT$) be represented
 by $X$ (resp. $T$).
Then \citesp{04SF}
\eqn{
\MMor_{\FFib/\cS}(\cT, \cX)/2\text{-isomorphisms}= \MMor_\cS(T, X).
}
Thus in this situation one may say that the only morphisms in
 $\cXT=\MMor_{\FFib/\cS}(\cT, \cX)$
 are identities.

Another way to explain this is that   when $X\in\cS$, $\cX=X\in \FFibS/\cS$,
 $\cX(\cT)=\MMor_{\FFib/\cS}(\cT, \cX)$ is a discrete category.

Then let $S=\Spec k$, $\cX=X$ a variety over $k$  and
 $\cA\ra \cS$ be $\Spec\bfA_k\ra \Spec k$ induced by the
 inclusion $k\subset \bfA_k$, and then we have $\cXA=\XA$ and
 $\cXS=\Xk$.
}

\defi{ \label{defi_obcat}
A category $\cXA'$ is an \emph{obstruction category} of $\cXA$, if it is
 a strictly full subcategory of $\cXA$ containing $\cXS$, in other words,
 if it  satisfies
\enmt{[\upshape (a)]
\item  $\cXS\subseteq \cXA'\subseteq \cXA$ and
\item \label{it_obcat_iso}
 if $x,y\in\cXA$ are isomorphic $($i.e., $2$-isomorphic in $\FFib/\cS$%
 $)$, then $x\in\cXA'$ if and only if $y\in\cXA'$.
}}

\defi{ \label{defi_pt_lift}
Let $\cX, \cY, \cT\in\FFib/\cS$ and  $f:\cY\ra \cX$ a $1$-morphism.
Let $\cYT'\subseteq \cYT$ be a full subcategory.
\enmt{[\upshape (a)]
\item Denote by $f(\cYT')$  the essential image of the functor
 $\MMor_{\FFib/\cS}(\cT, f): \cYT'\ra \cXT$.
\item We say that a point $x\in\cXT$ \emph{can be lifted} to a point in
  $\cYT'$ if $x\in f(\cYT')$.
 This is equivalent to say that
  there exists $y\in\cYT'$ such that $x$ and $f\circ y$ are
  $2$-isomorphic, that is, the following  diagram is $2$-commutative
 \eqn{\xymatrix{
  &\cT\ar[ld]_-{y} \ar[d]^-{x} \\
 \cY\ar[r]^-{f} &\cX
 }}
 In this case, $y$ is called a \emph{lift}
  of $x$ to $\cYT'$.
}}

The next lemma gives a condition describing whether a point $x\in\cXT$
 can be lifted to $\cYT$.
\lemm{ \label{lemm_sect_lift}
With notations in {\upshape \ref{defi_pt_lift}}, let $x\in\cXT$ be arbitrary
 and  $\cY_x=\cY\tm_{f, \cX, x} \cT$  the $2$-fibred product in $\FFib/\cS$.
Then $x\in f(\cYT)$ if and only if $\cY_x\ra \cT$ has a quasi-section in
 $\FFib/\cS$.
}
\pf{
By  properties of a $2$-fibred product, in the following $2$-commutative
 diagram
\eqn{\xymatrix{
\cT\ar@/^/[drr]^-\id \ar@{.>}[dr]|-s \ar@/_/[ddr]_-y \\
 &\cY_x\ar[r]\ar[d] &\cT\ar[d]^-x \\
 &\cY\ar[r]^-f &\cX
}}
 given any lift  $y\in \cYT$ of $x$, there is, up to a $2$-isomorphism, a unique
 quasi-section $s$ of $\cY_x\ra \cT$ making the whole diagram $2$-commutes.
Thus compositing with $\cY_x\ra\cY$ yields an one-to-one correspondence
\eqn{
\{\text{quasi-sections of }\cY_x\ra\cT\}/2\text{-isomorphisms}\os{\sim}{\lra}
 \{\text{lifts of $x$ to $\cYT$}\}/2\text{-isomorphisms}.
}
The proof is complete.
}

\subsection{The $F$-obstruction} \label{Fob}
Now we generalize the $F$-obstruction in \ref{Fob_c} defined on varieties to
 that on fibred categories.
Keep notations in \ref{X(A)_X(S)}.
\defi{
Let $F: \CC\ra\cD$ be a functor from a $2$-category to a category forgetting
 $2$-morphisms.
We say that $F$ is \emph{stable} if for any $1$-morphisms $f$ and $g$ of $\CC$
 that is $2$-isomorphic, we have $F(f)=F(g)$.
When $\cD$ is viewed as a $2$-category, this is equivalent to say that, $F$ is
 a $2$-functor.

Clearly if $F$ is stable, it remains so when restricted to a sub $2$-category
 of $\cC$.
}

\defi{ \label{defi_Fob}
Let $F: (\FFib/\cS)^\circ\ra \SSet$ be a stable functor.
Let $\cT$ be a fibred category over $\cS$ and $A\in F(\cX)$.
The \emph{evaluation} of a $\cT$-point
 $x\in \cXT=\MMor_{\FFib/\cS}(\cT, \cX)$
 at $A$ is defined to be  the image  of  $A$ under the pull-back map
 $F(x): F(\cX)\ra F(\cT)$ induced by $x$, denoted by $A(x)$.
}

Then forgetting morphisms in $\cXS$ and $\cXA$,
 we have a commutative diagram
\eqn{\xymatrix{
\cXS\ar[rr]^-{\MMor_{\FFib/\cS}(q, \cX)}\ar[d]^-{A(-)} &&\cXA\ar[d]^-{A(-)} \\
F(\cS)\ar[rr]^-{F(q)} &&F(\cA)
}}

\defi{
The \emph{obstruction given by $A$} is  the full
 subcategory $\cXA^A$ of $\cXA$ whose objects are characterized by
\eqn{
\cXA^A=\{x\in \cXA\mid A(x)\in \im F(q)\}.
}

The  \emph{$F$-category} is the full subcategory $\cXA^F$ whose objects are
 characterized by
\eqn{
\cXA^F=\bigcap_{A\in F(\cX)} \cXA^A =
 \{x\in \cXA\mid \im F(x)\subseteq \im F(q)\}.
}
 also denoted by $\cXA^{F(\cX)}$.
}

Then we have
\eqn{
\emptyset\subseteq \cXS\subseteq \cXA^F\subseteq \cXA^A\subseteq
 \cXA,
}
 and \ref{defi_obcat} \eqref{it_obcat_iso} holds.
In particular, both  $\cXA^F$ and $\cXA^A$ are obstruction categories.

\subsection{Cohomology for fibred categories} \label{coh}
For any  $\cT\in\FFib/\cS$,
 denote by $\cT_\et$ (resp. $\cT_\fppf$) the site with topology inherited from
 the big \'etale site $(\SSch/S)_\et$
 (resp. the big fppf site $(\SSch/S)_\fppf$) \citespp{067N}{06TP}.
Let  $\tau\in\{\et, \fppf\}$.
Then we have the cohomology functors
\eqn{
H^i(\cT_\tau, -): \AAb(\cT_\tau)\ra \AAb,
}
 where $\AAb(\cT_\tau)$ is the category of sheaves on $\cT_\tau$ and
 $\AAb$ the category of abelian groups.
See \citesp{075E}.
We will write $H_\tau^i(\cT, -)=H^i(\cT_\tau, -)$.
For $T\in\cS$, the corresponding $T\in\FFib/\cS$
 (see \ref{rk_emb} \eqref{it_U_Fib/S}) has a final object $T$.
Thus by  \citesp{075F}, $H_\tau^i(T, -)$  recovers the $\tau$-cohomology for
 schemes over $S$.
In particular,  we have  $H_\tau^i(\cS, -)=H_\tau^i(S, -)$, agrees with the
 $\tau$-cohomology for the scheme $S$.

\subsection{Functoriality of $\tau$-topoi} \label{pb_sh}
Let $f: {\cT'}\ra\cT$ be a $1$-morphism in $\FFib/\cS$.
Then compositing with $f$ gives the functor
 $f^*: \PPSh(\cT)\ra \PPSh(\cT')$, which has a right adjoint $f_*$.
Both $f^*$ and $f_*$ send $\tau$-sheaves to $\tau$-sheaves, and define
 a morphism of topoi
\eqn{
f=(f_*, f^*): \St{\cT'}\ra \St{\cT},
}
 which is also a pair of adjoint functors
 with $\SSh$ replaced by $\AAb$ (or $\SSGrp$ the
 category of sheaves of groups) and is compatible with the classical case when
 $\cT$ and $\cT'$ are represented by $S$-schemes.
See \citespp{06TS}{075J}.
A $2$-morphism $\xi: f\Ra g$ in $\FFib/\cS$
 canonically  induces a morphism of functors
\eq{ \label{eq_xi_upper_star}
\xi^*: g^*\Ra f^*,
}
 and hence a  morphism of functors $\xi_*: f_*\Ra g_*$ by adjunction.
In this way, $\xi$ induces a $2$-morphism of topoi
 $\xi=(\xi_*,\xi^*): f\Ra g$, which is actually a $2$-isomorphism.
The assignment is compatible with horizontal and vertical compositions of
 $2$-morphisms.
See \citespp{06TL}{06TM}.
\rk{
We prefer to use big $\tau$-topoi because they are suitable for our framework.
For example, the lisse-\'etale topos on algebraic stacks is not fuctorial
 \citesp{07BF} and the pull-back functor has to be  constructed using
 cohomological descent (c.f. \cite{olsson07sheaves, lo08six-i}).
Using big $\tau$ sites makes this simpler and works with arbitrary fibred
 categories.
See also Remark \ref{rk_do_bit_small}
}

\subsection{} \label{coh_func} Now using the notion of
 $F$-obstruction \ref{Fob},
 we generalize the cohomological obstructions in \ref{eg_cohob}.
Let $f: {\cT'}\ra\cT$ be a $1$-morphism in $\FFib/\cS$ and
 $\sF\in\AAb(\cT'_\tau)$.
Then we have an obvious isomorphism
\eqn{
\Gamma(\cT_\tau, f_*\sF)=\MMor_{\PPSh(\cT)}(e_\cT, f_*\sF)\cong
 \MMor_{\PPSh(\cT')}(f^*e_\cT, \sF)=
 \MMor_{\PPSh(\cT')}(e_{\cT'}, \sF)=\Gamma(\cT'_\tau, \sF),
}
 functorial in $\sF$, where  $e_\cT$ (resp. $e_{\cT'}$) is any final object in
 $\PPSh(\cT)$ (resp. $\PPSh(\cT')$).
Thus it induces the morphism of right derived functors
\eq{ \label{eq_pb_coh_der}
f^*: R\Gamma(\cT_\tau, -)\os{adj}{\Ra}
 R\Gamma(\cT_\tau, Rf_*f^*-)\os{\sim}{\Ra}
 R\Gamma(\cT'_\tau, f^*-): D(\cT_\tau)\ra D(\AAb),
}
 where $D(-)$ is the corresponding derived category.
For any $L\in D(\cT_\tau)$, we have a  canonical isomorphism
 $\Gamma(\cT_\tau, L)\os{\sim}{\ra} \Hom^._{K(\cT_\tau)}(\ZZ_\cT, L)$,
 where $K(-)$ is the corresponding homotopy category.
Passing to derived version, we have
 $R\Gamma(\cT_\tau, L)\os{\sim}{\ra} R\Hom^._{K(\cT_\tau)}(\ZZ_\cT, L)$.
Taking cohomology,  it becomes
 $H_\tau^i(\cT, L)\os{\sim}{\ra} \Hom_{D(\cT_\tau)}(\ZZ_\cT, L)$,
 which fits in to the following commutative diagram
\eq{ \label{eq_pb_coh_tps}
\xymatrix{
H_\tau^i(\cT, L)\ar[r]^-{\sim}\ar[d]^{f^*}
 &\Hom_{D(\cT_\tau)}(\ZZ_\cT, L)\ar[d]^-{f^*} \\
H_\tau^i(\cT', f^*L)\ar[r]^-{\sim}
 &\Hom_{D(\cT'_\tau)}(\ZZ_{\cT'}, f^*L)
}}
 functorial in $L$,  where
 the left $f^*$ is defined in \eqref{eq_pb_coh_der} and the right one
 is the functor $D(\cT_\tau)\ra D(\cT'_\tau)$ induced from
 $f^*: \St{\cT}\ra \St{\cT'}$.
In general, both $f^*$ in the diagram exists for $f$ being a
 $1$-morphism of topoi.
In particular, we obtain a morphism of functors
\eq{ \label{eq_pb_coh_gen}
f^*: H_\tau^i(\cT, -)\Ra H_\tau^i(\cT', f^*-): \AAb(\cT_\tau)\ra \AAb.
}
\lemm{ \label{lemm_coh_func}
With above notations, let
\eqn{\xymatrix@-2mm{
\cT'\rrtwocell^f_g{\xi} \ar[dr]_-{a'} & &\cT\ar[dl]^-a \\
 &\cS
}}
 be a $2$-isomorphism in $\FFib/\cS$.
\enmt{[\upshape (i)]
\item \label{it_2-mor_coh}
 We have the following commutative diagram of functors on $\AAb(\cT_\tau)$
 \eqn{\xymatrix@R=2mm{
  &H_\tau^i(\cT', g^*-) \ar@{=>}[dd]^-{\xi^*} \\
 H_\tau^i(\cT, -)\ar@{=>}[ru]^-{g^*} \ar@{=>}[rd]^-{f^*} \\
  &H_\tau^i(\cT', f^*-)
 }}
\item \label{it_coh_stable}
 For all $\sG\in\AAb(\cS_\tau)$, let $\sF=a^*\sG\in\AAb(\cT_\tau)$ be the
  pull-back sheaf. Then  we have $g^*\sF=f^*\sF\in\AAb(\cT'_\tau)$, and
  $\xi^*_\sF$ is the identity in \eqref{it_2-mor_coh}.
}}
\pf{
\eqref{it_2-mor_coh} comes from the interpretation of \eqref{eq_pb_coh_gen}
 by \eqref{eq_pb_coh_tps}  and the definition of $\xi^*$
 \eqref{eq_xi_upper_star}.

For \eqref{it_coh_stable}, note that $\xi^*_\sF$ is the composition of
 the identities $g^*\sF=g^*a^*\sG=a'^*\sG=f^*a^*\sG=f^*\sF$,
 and the result follows.
}

\subsection{Cohomological obstructions} \label{Hiob}
Now fix  $\sG\in\AAb(\cS_\tau)$ and let $\sG_\cT=a^*\sG$ be the  pull-back of
 $\sG$ to $\cT$.
Then  we have $\cG_{\cT'}=a'^*\sG=f^*a^*\sG=f^*\cG_\cT$,
 and the pull-back map \eqref{eq_pb_coh_gen} becomes
\eq{ \label{eq_pb_coh}
f^*: H_\tau^i(\cT, \sG_\cT)\ra H_\tau^i({\cT'}, \sG_{\cT'}),
}
 which makes $H_\tau^i(-, \sG_-): (\FFib/\cS)^\circ\ra \SSet$ a functor.
See also \ref{eq_pb_coh_tps}.
Since $\sG\in\AAb(\cS_\tau)$, by  \ref{lemm_coh_func} \eqref{it_coh_stable}
 we have $f^*=g^*: H_\tau^i(\cT, \sG_\cT)\ra H_\tau^i({\cT'}, \sG_{\cT'})$
 for any $2$-isomorphism $f\os{\sim}{\Ra} g$,
 which is to say that the functor
 $H_\tau^i(-, \sG_-): (\FFib/\cS)^\circ\ra \SSet$ is stable.
From now on, if there is  no other confusion, we simply write
 $\sG\in\AAb(\cT_\tau)$ for $\sG_\cT$.
Thus we can take $F=H_\tau^i(-, \sG)=H_\tau^i(-, \sG_-)$ in \ref{Fob} and
 obtain $\cXA^{H_\tau^i(\cX, \sG)}$, the $H_\tau^i(-, \sG)$-category.

\egu{[Brauer-Manin obstruction]
The  cohomological Brauer-Grothendieck group is also defined for any fibred
 category $\cT$ over $\cS$.
It is the \'etale cohomology group $\Br\cT:=H_\et^2(\cT, \bfG_{m, \cT})$, where
 $\bfG_{m, \cT}$ is the pull-back  sheaf of $\bfG_{m, S}$ to $\cT$.
When $\cT$ is a scheme, this agrees with  the cohomological
 Brauer-Grothendieck group.
The Brauer group of a stack is considered in, for example, Bertolin and
 Galluzzi \cite[Thm. 3.4]{bg19gerbes}, and Zahnd \cite[4.3]{zahnd03gerbe-br}.

Let $F=\Br=H_\et^2(-, \bfG_m)$,  we obtain the \emph{Brauer-Manin category}
 $\cXA^{\Br}$, which coincides  with the usual Brauer-Manin set in the
 classical case.
}

\section{The generalized descent obstruction and descent by torsors}
 \label{desc}
In order to generalize the classical descent theory from varieties to fibred
 categories, we need to reformulate everything in details.
Since groups are not necessarily commutative, we can not use cohomology
 functors on abelian sheaves discussed in \ref{coh}, and shall do in a
 different manner.
Let us first recall the
\defi{[Torsors] \label{defi_torsor}
Let $\cC$ be a site and $\SSh(\cC)$ the correspondence topos.
\enmt{[\upshape (a)]
\item Let $\sU\in\SSh(\cC)$. A \emph{$\sU$-group}  is a group object in  the
  localized category $\SSh(\cC)/\sU$.
 Let $\sG$ be a $\sU$-group. An object $\sY\in\SSh(\cC)/\sU$  is a
  (right) \emph{$\sG$-object} over $\sU$ if it is acted by $\sG$ on the right.
\item Let $\sY$ be a $\sG$-object over $\sU$. Then $\sY$ is a (right)
 \emph{$\sG$-torsor over $\sU$} if the following two conditions hold
  in $\SSh(\cC)$:
 \enmt{[\upshape (i)]
 \item the morphism $\sY\ra\sU$ is an epimorphism, and
 \item the morphism
  \aln{
  \sY\tm_\sU\sG &\ra \sY\tm_\sU\sY \\
  (y, g) &\mpt (y, yg)
  }
  is an isomorphism.
 }
 Morphisms of $\sG$-torsors over $\sU$ are morphisms in $\SSh(\cC)/\sU$
  compatible with the action of $\sG$.
 Let $\sY$ be a $\sG$-torsor over $\sU$.
 If $\sY\ra \sU$ has a section in $\SSh(\cC)$, we call $\sY$ \emph{trivial},
  which is
  equivalent to say that $\sY\cong\sG_r$ as $\sG$-torsor over $\sU$,
  where $\sG_r$ is the torsor $\sG$
  acted on by itself on the right.
 Left torsors, left objects and $\sG_l$ are are defined in a similar way.
\item For $U\in\cC$,  we write $\ul U=h_U^\#\in\SSh(\cC)$
  for the sheafification of
  the  presheaf  $h_U=\MMor\cC(-, U)\in\PPSh(\cC)$.
 Then for any $\sG\in\SSGrp(\cC)$,
  by a \emph{$\sG$-torsor over $U$} we mean a $\sG$-torsor over
  $\ul U$.
 Denote by $\Tors(\cC/U, \sG)$ the category of $\sG$-torsors over $U$, which
  is actually a groupoid.
\item  Let   $e\in\SSh(\cC)$ be a final object of the topos $\SSh(\cC)$ and
  $\sG$ a $e$-group, which is to say that $\sG$ is an object of $\SSGrp(\cC)$.
 A \emph{$\sG$-torsor over $\cC$} is a $\sG$-torsor over $e$.
 In particular, it is an object in $\SSh(\cC)$.
 Denote by $\Tors(\cC, \sG)$ the category of $\sG$-torsors over $\cC$.
\item Denote by $\sfH^1(\cC/U, \sG)$ (resp. $\sfH^1(\cC, \sG)$) the set of
  isomorphism classes in $\Tors(\cC/U, \sG)$ (resp. $\Tors(\cC, \sG)$).
 All trivial torsors form an isomorphism class,
  called the \emph{neutral} element, which
  makes $\sfH^1(\cC/U, \sG)$ (resp. $\sfH^1(\cC, \sG)$) a  pointed set.
}}

\subsection{Torsors as fibred categories} \label{tors_as_fib}
Let  $\tau\in\{\et, \fppf\}$ and  $\cT\in\FFib/\cS$.
Recall that $\cT_\tau$ is with the inherited topology from
 $\cS_\tau=(\SSch/S)_\tau$.
Let  $\sG\in\SSGrp(\cT_\tau)$  and $\sY\in\Tors(\cT_\tau, \sG)$.
As presheaf, $\sY$ is embedded in  $\FFib/\cT$ by the functor
 $\eta_\cT$ in \eqref{eq_emb}.
Denote by $f: \cY\os{\sG}{\ra} \cT_\tau$ the image $\eta_\cT(\sY)\in\FFib/\cT$,
 called the fibred category \emph{associated to}  $\sY$,
 and we simply write  $f: \cY\ra \cT\in \Tors(\cT_\tau, \sG)$ to indicate it.
By  \ref{eq_emb},
 if $\sY\os{\sim}{\ra}\sZ$ is an isomorphism in $\Tors(\cT_\tau, \sG)$,
 then we have the corresponding isomorphism $\cY\os{\sim}{\ra}\cZ$
 in $\FFib/\cT$.
Keep in mind that $(\cY\ra\cT\ra\cS)\in\FFib/\cS$ and the topology inherited
 from  $\cT_\tau$ agrees with $\cY_\tau$ (inherited directly from $\cS_\tau$).

\subsection{Classify torsors} \label{tors_H1}
With notations in \ref{tors_as_fib}, let $T\in\cT$.
Then the first (non-abelian) \v Cech cohomology $\check H_\tau^1(\cT/T, \sG)$
 is defined to be  $\varinjlim_{\cU} \check H_\tau^1(\cU/T, \sG)$
 (where $\cU$ runs over all covering of $T$ in $\cT_\tau$) and
 there is a one-to-one correspondence
 (c.f. Giraud \cite[III.3.6.5 (5)]{giraud71cohnonab}) of pointed sets
\eq{ \label{eq_tors_H1}
\sfH^1(\cT_\tau/T, \sG) \os{\sim}{\ra} \check H_\tau^1(\cT/T, \sG).
}
If $\cT$ has a final object, then $\check H_\tau^1(\cT/T, -)$, up to an
 isomorphism, is independent of the choice of the final object $T\in\cT$.
This suggests us to write
 $\check H_\tau^1(\cT, -)$ for $\check H_\tau^1(\cT/T, -)$.
Note that in this case, $\sfH^1(\cT_\tau, -)=\sfH^1(\cT_\tau/T,-)$
 by definition.
It follows that \eqref{eq_tors_H1} becomes
\eq{ \label{eq_tors_H1_final}
\sfH^1(\cT_\tau, \sG) \os{\sim}{\ra} \check H_\tau^1(\cT, \sG),
}
 functorial in $\sG$.

In the case of $\sG\in\AAb(\cT_\tau)$, dropping the assumption that $\cT$ has
 a final object, we always  have an bijection (c.f.
 \cite[IV.3.4.2 (i)]{giraud71cohnonab})
\eq{ \label{eq_tors_H1_ab}
\sfH^1(\cT_\tau, \sG) \os{\sim}{\ra} H_\tau^1(\cT, \sG),
}
 which maps the neutral element to $0$ and is functorial in $\sG$.

\subsection{The pull-back of a torsor} \label{pb_tors}
From now on, we shall write $\sfH_\tau^1(\cT, -)$ for $\sfH^1(\cT_\tau, -)$.
Let $\sG\in\SSGrp(\cT_\tau)$  and  $\sY\in\Tors(\cT_\tau, \sG)$.
For any $1$-morphism $u: {\cT'}\ra\cT$, the pull-back sheaf
 $u^*\sY\in\St{\cT'}$ is a $u^*\sG$-torsor over $\cT'$, called the
 \emph{pull-back} of $\sY$ under $u$.
Thus we have a well-defined map
\al{
u^*: \sfH_\tau^1(\cT, \sG) &\ra \sfH_\tau^1(\cT', u^*\sG)
 \label{eq_pb_tors_gen} \\
 [\sY] &\mpt [u^*\sY]. \nonumber
}

The pull-back of a torsor by a $1$-morphism in $\FFibP/\cS$ commutes with the
 formation of the fibred category associated to it.
To be precise, we have the following
\lemm{ \label{lemm_pb_eta}
With notations in {\upshape \ref{pb_tors}}, suppose that $\cT',
 \cT\in\FFibP/\cS$.
Then the pull-back functor $u^*$ commutes with $\eta$ in \eqref{eq_emb}.
In other words, there is an equivalence in $\FFib/\cT'$
\eqn{
\alpha: \eta_{\cT'}(u^*\sY)\os{\approx}{\ra} u^*\eta_\cT(\sY),
}
 where the second $u^*$ is defined by {\upshape \ref{defi_pb_fib}}.
}
\pf{
Let $\cY=\eta_\cT(\sY)=\cT_\sY$, $\cY'=u^*\cY$ and
 $\cY''=\eta_{\cT'}(u^*\sY)=\cT'_{u^*\sY}$.
By the constructions of $\cT_\sY$ and
 $\cT'_{u^*\sY}$ in \ref{lemm_emb} \eqref{it_eta},
 and the definition of $u^*\sY$ (see \ref{pb_sh}), we may define  a
 $1$-morphism $u'': \cY''\ra \cY$ in $\FFib/\cS$,
 such that the outer square of the following diagram in $\FFib/\cS$
 is $2$-commutative
\eqn{\xymatrix{
\cY''\ar@/^/[drr]^-{f''} \ar@{.>}[dr]|-{\alpha} \ar@/_/[ddr]_-{u''} \\
 &\cY'=\cY\tm_\cT \cT'\ar[r]_-{f'}\ar[d]^-{u'} &\cT'\ar[d]^-u \\
 &\cY\ar[r]^-f &\cT
}}
 where $f'': \cY''\ra \cT'$ is  the structure morphism.
By the universal  property of $2$-fibred product, there exists a $1$-morphism
 $\alpha: \cY''\ra \cY'$ such that the whole diagram is $2$-commutative.
Note that $\alpha$ is actually a $1$-morphism in $\FFib/\cT'$ since
 $\MMor_{\FFib/\cS}(\cY'', \cT')$ is a discrete category.
By \citesp{003Z}, it suffices to show that
 $\alpha_{t'}: \cY''_{t'}\ra \cY'_{t'}$ is an  equivalence for any $t'\in\cT'$.

By \ref{lemm_pb_fib}, for  $t'\in\cT'$ and $t=u(t')\in\cT$, we have
 $\cY'_{t'}\os{\approx}{\ra} \cY_t$ induced by $u'$.
On the other hand, $\cY_t\cong sY(t)$ and $\cY''_{t'}\cong u^*\sY)(t')=\sY(t)$
 by the definition of $\eta$ and $u^*$.
Thus $u''$ induces an equivalence $\cY''_{t'}\os{\approx}{\ra} \cY_t$.
It follows that
 $\cY''_{t'}\os{\approx}{\ra} \cY'_{t'}$ for any $t'\in\cT'$,
 and one  checks that up to equivalences, this is exactly $\alpha_{t'}$.
The proof is complete.
}
\rk{
From the proof one sees that the result holds for arbitrary category $\cS$ and
 $\sY\in\PPSh(\cT)$.
}
\cor{ \label{cor_comp_fibp}
If $\cT\ra\cS\in\FFibP/\cS$ and $\cY\ra\cT\in\FFibP/\cT$, then
 $(\cY\ra\cT\ra\cS)\in\FFibP/\cS$.
}
\pf{
Since  $\cT\ra\cS\in\FFibP/\cS$, by definition, there exist
 $\cT'\ra\cS\in\FFibS/\cS$ and  an equivalence
 $\cT'\os{\approx}{\ra}\cT$ in  $\FFib/\cS$.
Then Lemma \ref{lemm_pb_eta} tells us that the $2$-base change
 $\cY'=\cY\tm_\cT \cT'\ra \cT'\in\FFibP/\cT'$, and
 $\cY'\os{\approx}{\ra}\cY$ in  $\FFib/\cS$.
Clearly we have $(\cY'\ra\cT'\ra\cS)\in\FFibP/\cS$.
It follows that $(\cY\ra\cT\ra\cS)\in\FFibP/\cS$.
}

Now we assume that  $\sG\in\SSGrp(\cS_\tau)$ (whose pull-back to any
 $\cT_\tau\in\FFib/\cS$ is also denoted by $\sG$), then \eqref{eq_pb_tors_gen}
 becomes
\eq{ \label{eq_pb_tors}
u^*: \sfH_\tau^1(\cT, \sG) \ra \sfH_\tau^1(\cT', \sG),
}
 which makes $\sfH_\tau^1(-, \sG): (\FFib/\cS)^\circ\ra \SSet$ a functor.
Then \eqref{eq_tors_H1_final} (resp. \eqref{eq_tors_H1_ab}) is also functorial
 in $\cT\in\FFibR/\cS$ (resp. $\FFib/\cS$ if $\sG\in\AAb(\cS_\tau)$).

\subsection{} Let $\sG\in\SSGrp(\cS_\tau)$.
Since in general, $\cT\in\FFib/\cS$ need not to have a final object, we can not
 choose $F$ in \ref{Fob} to be the first \v Cech cohomology  in \ref{tors_H1}.
Instead we shall use $\sfH_\tau^1(-, \sG)$ directly.
First we need to verify that
\lemm{ \label{lemm_sfH1_stable}
If  $\sG\in\SSGrp(\cS_\tau)$, then
 the functor  $\sfH_\tau^1(-, \sG):  (\FFib/\cS)^\circ\ra \SSet$  is stable.
}
\pf{
Let
\eqn{\xymatrix@-2mm{
\cT'\rrtwocell^f_g{\xi} \ar[dr]_-{a'} & &\cT\ar[dl]^-a \\
 &\cS
}}
 be a $2$-isomorphism in $\FFib/\cS$.
Then we have the following commutative diagram of functors on
 $\SSGrp(\cT_\tau)$
\eqn{\xymatrix@R=2mm{
 &\sfH_\tau^1(\cT', g^*-) \ar@{=>}[dd]^-{\xi^*} && &[g^*\sY]\ar@{|->}[dd] \\
\sfH_\tau^1(\cT, -)\ar@{=>}[ru]^-{g^*} \ar@{=>}[rd]^-{f^*} &&
 &[\sY]\ar@{|->}[ru] \ar@{|->}[rd] \\
 &\sfH_\tau^1(\cT', f^*-) && &[f^*\sY]
}}
Note that $\sG\in\SSGrp(\cS_\tau)$.
Then $g^*\sG=f^*\sG\in\SSGrp(\cT'_\tau)$ and
 $\xi^*_\sY: g^*\sY\os{\sim}{\ra} f^*\sY$ in $\SSh(\cT'_\tau)$.
If follows that $[g^*\sY]=[f^*\sY]$ in  $\sfH_\tau^1(\cT',\sG)$.
Thus $\xi^*$ in the above diagram is the identity.
The proof is complete.
}

Then we have  $\cXA^{\sfH_\tau^1(\cX, \sG)}$,
 the $\sfH_\tau^1(-, \sG)$-category, and the following
\defi{ \label{defi_desc}
Define the \emph{descent category} to be the full
 subcategory $\cXA^\desc$ of $\cXA$ whose objects are characterized by
\eqn{
\cXA^\desc=\bigcap_{\sG\in\SSGrp(\cS_\fppf)} \cXA^{\sfH_\fppf^1(\cX, \sG)}.
}}
Clearly, $\cXA^\desc$ is an obstruction category.

\subsection{} We shall first compare $\desc$ with $\cdesc$ \eqref{eq_cdesc}.
For  $\cT\in\FFib/\cS$ represented by $T$, $T$ is a final object of $\cT$.
Recall in \ref{tors_H1} that in this case
 we write $\check H_\tau^1(\cT, \sG)=\check H_\tau^1(\cT/T, \sG)$.
Suppose that $\sG\in\SSGrp(\cS_\tau)$.
Then $\check H_\tau^1(-, \sG)$ is identified with
 $\sfH_\tau^1(-,\sG)$ via \eqref{eq_tors_H1_final},
 and thus by   \ref{lemm_sfH1_stable} it is also
 a stable functor on $\FFibR/\cS$ (see \ref{lemm_emb} \eqref{it_fib_rep}).
Thus $\cXA^{\check H_\tau^1(\cX, \sG)}$ is defined if $\cX,\cA\in\FFibR/\cS$.
\prop{ \label{prop_cXAH1}
With notations above, let $\sG\in\SSGrp(\cS_\tau)$.
\enmt{[\upshape (i)]
\item \label{it_cXAH1_rep}
 If $\cX,\cA\in\FFibR/\cS$,
  then $\cXA^{\sfH_\tau^1(\cX, \sG)}=\cXA^{\check H_\tau^1(\cX, \sG)}$.

 In particular, in the classical case where  $X$  is a variety over a number
  field $k$,
  the set $\XA^\desc$ defined by {\upshape \ref{defi_desc}} is contained in the
  classical  descent set defined by {\upshape \eqref{eq_cdesc}}, i.e.
 \eq{ \label{eq_desc_in_cdesc}
 \XA^\desc\subseteq \XA^\cdesc.  
 }
\item \label{it_cXAH1_ab}
 If  $\sG\in\AAb(\cS_\tau)$,
  then $\cXA^{\sfH_\tau^1(\cX, \sG)}=\cXA^{H_\tau^1(\cX, \sG)}$.
}}
\pf{
In \eqref{it_cXAH1_rep}, all of $\cX$, $\cA$ and $\cS$ are in $\FFibR/\cS$.
Since $\check H_\tau^1(-, \sG)$ is identified with $\sfH_\tau^1(-, \sG)$ via
 \eqref{eq_tors_H1_final}, the resulting obstruction categories are the same.
This shows the first equality.

In the classical case,  let $S=\Spec k$, and
 $\cA=\Spec \bfA_k, \cX=X\in\FFibR/\cS$.
Note that in \eqref{eq_cdesc}, $G$ runs over all affine $k$-groups, viewed as a
 sheaf $\sG\in\SSGrp(\cS_\fppf)$.
Thus the inclusion \eqref{eq_desc_in_cdesc} follows.

The statement \eqref{it_cXAH1_ab} follows from the isomorphism
 \eqref{eq_tors_H1_ab} functorial in $\cT\in\FFib/\cS$.
}

\subsection{The obstruction given by torsors} \label{ob_tors}
We now extend the classical theory of descent by torsors.
Recall that a torsor $\sY\in\Tors(\cT_\tau, \sG)$ is mapped by the  functor
 $\eta_\cT$ to  $f: \cY\ra \cT\in \FFib/\cS$ associated to it
 (see \ref{tors_as_fib}), and represents a class
 $A=[\sY]\in\sfH_\tau^1(\cT, \sG)$.
We write $\cXA^f=\cXA^A$, the obstruction category given by $A$.
\prop{ \label{prop_dbt}
Suppose that $\cA\in\FFibS/\cS$ and $\cX\in\FFibP/\cS$.
Let $\sG\in\SSGrp(\cS_\tau)$ and  $f: \cY\ra \cX\in \Tors(\cX_\tau, \sG)$.
Then $\Ob(\cXA^f)$ is characterized by
\eqn{
\cXA^f=\bigcup_{\s\in \check H_\tau^1(\cS, \sG)} f^\s(\cY^\s(\cA)),
}
 where $f^\s: \cY^\s\ra \cX\in\Tors(\cX_\tau, \sG^\s)$ is the twist of
 $\cY$ by $\s$ {\upshape (}see {\upshape \ref{defi_twist_tors}\upshape )}.
}
This leads to the following
\cor{ \label{cor_desc_by_dbt}
Suppose that $\cA\in\FFibS/\cS$ and $\cX\in\FFibP/\cS$.
Then we have
\eqn{
\cXA^\desc=
 \bigcap_{\os{\sG\in\SSGrp(\cS_\fppf)}{f: \cY\ra \cX\in\Tors(\cX_\fppf, \sG)}}
 \bigcup_{\s\in \check H_\fppf^1(\cS, \sG)} f^\s(\cY^\s(\cA)).
}}
\pf{
Since $\sfH_\fppf^1(\cX, \sG)$ is the pointed set of isomorphism classes of
 $\Tors(\cX_\fppf, \sG)$ (see \ref{ob_tors}), the result
 follows immediately from \ref{prop_dbt}.
}

\cor{ \label{cor_rational_by_dbt}
With assumptions and notations in {\upshape \ref{prop_dbt}}, we have
\eqn{
\cXS=\bigcup_{\s\in \check H_\tau^1(\cS, \sG)} f^\s(\cY^\s(\cS)).
}}
\pf{
Take $\cA=\cS$ in \ref{prop_dbt}, and note that $\cXS^f=\cXS$.
Then the result follows.
}

The proof of \ref{prop_dbt} will be given in \ref{pf_prop_dbt}.
Note that the assumption that
 $\cA\in\FFibS/\cS$ and $\cX\in\FFibP/\cS$
 holds automatically in the
 classical  case where $X$ is a $k$-variety.
See \ref{lemm_emb}.

\subsection{Bitorsors and the inverse torsor} \label{bitors_invtors}
Let $\cT\in\FFib/\cS$ and $\sG,\sH\in\SSGrp(\cT_\tau)$.
For $\sY\in\Tors(\cT_\tau, \sG)$,
 if $\sH$ acts on $\sY$ on the left such that $\sY$ is also a left
 $\sH$-torsor over $\cT_\tau$, $\sY$ is called a \emph{$\sH$-$\sG$-bitorsor}.

Now suppose that $\sY$ is a $\sH$-$\sG$-bitorsor.
Let  $\sG$ and $\sH$ act on $\sY$ by
\aln{
\sG\tm \sY\tm \sH &\ra \sY \\
(g, y, h) &\mpt h^{-1}yg^{-1}.
}
Then $\sY$ becomes a $\sG$-$\sH$-bitorsor, called the \emph{inverse torsor}
 of $\sY$, denoted by $\sY^\circ$.
Of cause, we have $\sY^\circ\in\Tors(\cT_\tau, \sH)$.

In particular, for arbitrary $\sY\in\Tors(\cT_\tau, \sG)$,
 let $\ad(\sY)=\Aut_\sG(\sY)$ be the adjoint group of $\sY$,
 so that it acts $\sY$ on the left.
Then $\sY$ becomes an  $\ad(\sY)$-$\sG$-bitorsor, so that $\sY^\circ\in
 \Tors(\cT_\tau, \ad(\sY))$.
See \cite[III.1.5]{giraud71cohnonab} for more details.

\subsection{The contracted product of torsors} \label{contr_tors}
Let $\sF,\sG,\sH\in\SSGrp(\cT_\tau)$.
Let  $\sY$ (resp. $\sZ$) be a $\sF$-$\sG$- (resp. $\sG$-$\sH$-) bitorsor.
The \emph{contracted product} \cite[III.1.3]{giraud71cohnonab}
 of $\sY$ and $\sZ$ is the quotient of $\sY\tm\sZ$ by
 the diagonal operation of $\sG$
\aln{
\sY\tm\sZ\tm\sG &\ra \sY\tm\sZ \\
(y, z, g) &\mpt (yg, g^{-1}z),
}
 which  is a $\sF$-$\sH$-bitorsor over $\cT_\tau$,
 denoted by $\sY\os{\sG}{\tm}\sZ$.
The contracted product is associative, i.e., if given $\sE\in\SSGrp(\cT_\tau)$
 and an $\sE$-$\sF$-bitorsor $\sX$, we have an canonical isomorphism of
 $\sE$-$\sH$-bitorsors \cite[III.1.3.5]{giraud71cohnonab}
\eq{ \label{eq_contr_asso_tors}
(\sX\os{\sF}{\tm}\sY)\os{\sG}{\tm}\sZ\os{\sim}{\ra}
 \sX\os{\sF}{\tm}(\sY\os{\sG}{\tm}\sZ).
}
Also it commutes with the pull-back, that is,
 for $1$-morphism $u: \cT'\ra \cT$, we have an canonical isomorphism of
 $u^*\sF$-$u^*\sH$-bitorsors
\eq{ \label{eq_contr_pb_tors}
u^*(\sY\os{\sG}{\tm}\sZ)\os{\sim}{\ra} u^*\sY\os{u^*\sG}{\tm}u^*\sZ.
}

\lemm{ \label{lemm_contr_inv}
Let $\sY$ be a $\sH$-$\sG$-bitorsor.
Then we have isomorphisms
\gan{
\sY\os{\sG}{\tm}\sG_b\cong\sY, \quad  \sH_b\os{\sH}{\tm}\sY\cong\sY \quad
 (\text{as $\sH$-$\sG$-bitorsors}), \\
\sY\os{\sG}{\tm}\sY^\circ\cong\sH_b \quad(\text{as $\sH$-bitorsors}), \\
\sY^\circ\os{\sH}{\tm}\sY\cong\sG_b \quad(\text{as $\sG$-bitorsors}),
}
 where $\sH_b$ and $\sG_b$ are the bitorsors by the obvious actions.
}
\pf{
See \cite[III.1.3.1.3 and III.1.6.5]{giraud71cohnonab}.
}

\subsection{} \label{twist_tors}
Let
\eqn{\xymatrix@-2mm{
\cT'\ar[rr]^g \ar[dr]_-{a'} & &\cT\ar[dl]^-a \\
 &\cS
}}
 be a $1$-morphism in $\FFib/\cS$.
Then we have
\eq{ \label{eq_pb_comp_H1}
g^*\circ a^*=a'^*
}
 for the  pull-back maps defined by \ref{pb_tors}.

\defi{[Twist of a torsor] \label{defi_twist_tors}
With notations in \ref{twist_tors},
 let $\sG\in\SSGrp(\cS_\tau)$ (whose pull-back to any $\cT_\tau\in\FFib/\cS$
 is also denoted by $\sG$),
 $\s\in \check H_\tau^1(\cS, \sG) \os{\sim}{\ra} \sfH_\tau^1(\cS, \sG)$
 (see \eqref{eq_tors_H1_final}),
 and $\sZ\in\Tors(\cS_\tau, \sG)$ represents  $\s$.
For any  $\sY\in\Tors(\cT_\tau, \sG)$ with $\eta_\cT(\sY)= (f: \cY\ra \cT)$
 associated to it,
 we know the class
 $[\sY\os{\sG}{\tm}(a^*\sZ)^\circ]$ is independent of the choice of $\sZ$
 (see \ref{contr_tors}).
A \emph{twist} of $\sY$ by $\s\in \check H_\tau^1(\cS, \sG)$ is any
 representative of the class $[\sY\os{\sG}{\tm}(a^*\sZ)^\circ]$,
 denoted by $\sY^\s\in\Tors(\cT_\tau, \sG_\cT^\s)$,
 where $\sG_\cT^\s\in\SSGrp(\cT_\tau)$, depending on $\sZ$, is any
 representative of the
 isomorphism class $[a^*\ad(\sZ)]$ (also independent of $\sZ$) of groups.
We also write $f^\s: \cY^\s\ra\cT$ for  $\eta_\cT(\sY^\s)$ associated to
 $\sY^\s$ (see \ref{tors_as_fib}).
Up to isomorphisms  in $\FFib/\cT$,  $f^\s: \cY^\s\ra\cT$ is determined only
 by $\sY$ and $\s$.
}

\lemm{ \label{lemm_pb_twist_tors}
With the above notations, we have $g^*(\sG_\cT^\s)\cong \sG_{\cT'}^\s$ in
 $\SSGrp(\cT'_\tau)$
 $($this suggests us to write $\sG^\s$ for $\sG_\cT^\s$ for
 any $\cT\in\FFib/\cS$ if no confusion caused$)$
 and $(g^*\sY)^\s\cong g^*(\sY^\s)$ in $\Tors(\cT'_\tau, \sG^\s)$.
}
\pf{
Let $\sZ\in\Tors(\cS_\tau, \sG)$ represents $\s\in\sfH_\tau^1(\cS, \sG)$.
Then by \eqref{eq_pb_comp_H1} and \ref{defi_twist_tors},
\eqn{
\sG_{\cT'}^\s\cong a'^*\ad(\sZ)= (g^*\circ a^*)\ad(\sZ)\cong g^*\sG_\cT^\s
}
 and the first statement follows.

For the second  statement,
 by \eqref{eq_contr_pb_tors}, \eqref{eq_pb_comp_H1}
 and \ref{defi_twist_tors}, we have
\eqn{
(g^*\sY)^\s=g^*\sY\os{\sG}{\tm}(a'^*\sZ)^\circ=
 g^*\sY\os{\sG}{\tm} g^*(a^*\sZ)^\circ\cong
 g^*(\sY\os{\sG}{\tm}(a^*\sZ)^\circ)=g^*(\sY^\s),
}
 compatible with the action of $\sG^\s$.
The proof is complete.
}

\subsection{Proof of $\ref{prop_dbt}.$ Proposition} \label{pf_prop_dbt}
Let $\sY\in\Tors(\cX_\tau, \sG)$ whose image under $\eta_\cX$ is
 $f: \cY\ra \cX$.
By definition, $x\in \cXA^f$ means that the evaluation of
 $[\sY]\in \sfH_\tau^1(\cX, \sG)$ at $x$ comes from $\sfH_\tau^1(\cS,\sG)$,
 that is, there exists
 $\s\in \check H_\tau^1(\cS, \sG) \os{\sim}{\ra} \sfH_\tau^1(\cS, \sG)$ such
 that
\eq{ \label{eq_Y(x)=s}
[\sY](x)=q^*\s \in \sfH_\tau^1(\cA, \sG).
}
Let $\sZ\in\Tors(\cS_\tau, \sG)$ represents $\s$.
Then \eqref{eq_Y(x)=s} means that $x^*\sY\cong q^*\sZ$ in
 $\Tors(\cA_\tau, \sG)$,
 which, by \ref{defi_twist_tors} and \ref{lemm_contr_inv}, implies that
\eqn{
(x^*\sY)^\s:=x^*\sY\os{\sG}{\tm}(q^*\sZ)^\circ\cong
 q^*\sZ\os{\sG}{\tm}(q^*\sZ)^\circ\cong(\sG^\s)_r
}
 in $\Tors(\cA_\tau, \sG^\s)$ (the notation for $\sG^\s$ is in
 \ref{lemm_pb_twist_tors}).
Conversely, if  $x^*\sY\os{\sG}{\tm}(q^*\sZ)^\circ\cong(\sG^\s)_r$
 in $\Tors(\cA_\tau, \sG^\s)$
 then by \ref{lemm_contr_inv} again and \eqref{eq_contr_asso_tors},
\eqn{
x^*\sY\cong x^*\sY\os{\sG}{\tm}\sG_l\cong x^*\sY\os{\sG}{\tm}(q^*\sZ)^\circ
 \os{\sG^\s}{\tm}q^*\sZ \cong (\sG^\circ)_r\os{\sG^\s}{\tm}q^*\sZ\cong
 q^*\sZ
}
 in  $\Tors(\cA_\tau, \sG)$.
On the other hand, by \ref{lemm_pb_twist_tors},
 $x^*(\sY^\s)\cong (x^*\sY)^\s$.
It follows that  \eqref{eq_Y(x)=s} is equivalent to that
 $x^*(\sY^\s)\in\Tors(\cA_\tau, \sG^\s)$ is trivial, i.e., $x^*(\sY^\s)\ra e$
 has a section, where $e$ is a final object of $\St{\cA}$.
Note that  $\cX, \cA\in\FFibP/\cS$, by \ref{lemm_pb_eta},
 there is an equivalence in $\FFib/\cA$
\eq{ \label{eq_pb_eta}
\eta_\cA(x^*(\sY^\s))\os{\approx}{\ra}(\cY^\s)_x,
}
 where $(\cY^\s)_x:= x^*(\cY^\s)\in\FFib/\cA$.
Thus that $x^*(\sY^\s)\ra e$ has a section
\aln{
&\Leftrightarrow \eta_\cA(x^*(\sY^\s))\ra \cA \text{ has a section in }
 \FFib/\cA &(\text{$\eta_\cA$ is an equivalence}), \\
&\Leftrightarrow (\cY^\s)_x\ra \cA \text{ has a section in }
 \FFib/\cA &(\text{by \eqref{eq_pb_eta}}), \\
&\Leftrightarrow (\cY^\s)_x\ra \cA \text{ has a quasi-section in }
 \FFib/\cS &(\text{by $\cA\in\FFibS/\cS$ and \ref{rk_emb} \eqref{it_qs_s}}), \\
&\Leftrightarrow x\in f^\s(\cY^\s(\cA))
 &(\text{by \ref{lemm_sect_lift}}).
}
The proof is complete. \qed

\section{The second descent obstruction and descent by gerbes} \label{sdesc}
We now turn back to abelian sheaves.
Under the discussion in \ref{Fob} and \ref{coh}, we define the second
 descent obstruction on a fibred category and formulate the idea of
 descent by gerbes.

\subsection{The obstruction given by second cohomology} \label{ob_scoh}
Let $\tau\in\{\et,\fppf\}$, $\sG\in\AAb(\cS_\tau)$ (whose pull-back to any
 $\cT_\tau$ is also denoted by $\sG$) and $F=H_\tau^2(-, \sG)$, as in
 \ref{coh}.
We have $\cXA^{H_\tau^2(\cX, \sG)}$, the $H_\tau^2(-, \sG)$-category.
Elements of $H_\tau^2(\cX, \sG)$ are in one-to-one correspondence
 with those of $\sfH_\tau^2(\cX, \sG)$,
 the $\sG$-equivalence classes of gerbes over $\cX_\tau$ bounded by $\sG$
 (see \ref{defi_gerbe} and  \ref{gerb_H2} below).
Thus any $A\in H_\tau^2(\cX, \sG)$ is represented by a gerbe
 $f: \cY\ra \cX\in\Gerb(\cX_\tau, \sG)$.
That is,  the class $[\cY]$ maps to $A$.
See \eqref{eq_gerb_H2}.
Then we  have $\cXA^f:=\cXA^A$, the obstruction given by $A$.

\defi{ \label{defi_sdesc}
Define the \emph{second descent category} to be the full
 subcategory $\cXA^\sdesc$ of $\cXA$ whose objects are characterized by
\eqn{
\cXA^\sdesc=\bigcap_{\sG\in\AAb(\cS_\et)} \cXA^{H_\et^2(\cX, \sG)}.
}}
It is clear that $\cXA^\sdesc$ is an obstruction category satisfying  that
 $\cXA^\sdesc\subseteq \cXA^{H_\et^2(\cX,\sG)}$ for all $\sG\in\AAb(\cS_\et)$.
In particular, we have $\cXA^\sdesc\subseteq \cXA^{\Br}$ since $\bfG_{m, S}$
 represents a sheaf in $\AAb(\cS_\et)$.
Unfortunately,
 we do not know the relation between $\cXA^\sdesc$ and $\cXA^\desc$.
In next two sections, we shall obtain some full subcategories between $\cXS$
 and $\cXA^\desc$.
\rk{
We prefer to use \'etale topology for $\sdesc$ since it is easy to compare with
 $\Br$.
In fact we may also define $\sdesc$ replacing $\et$ by $\fppf$ and most of
 the results  involving $\sdesc$ in this text remain correct.
}

We now recall the
\defi{ \label{defi_gerbe}
Let $\cC$ be a site.
\enmt{[\upshape (a)]
\item \label{it_defi_stack}
 A \emph{stack in groupoid} over $\cC$ is a  fibred category
  $\cY\ra \cC$ such that for any covering $\cU$ in $\cC$, the functor
  $\cY_\cU\ra DD(\cU)$ which associates to an object its \emph{canonical
  descent datum}, is an equivalence \citesp{026E}.
 In this text, by a stack we  always mean a stack in groupoid.
 Denote by $\SStk/\cC$ (resp. $\SStkS/\cC$, resp. $\SStkP$) the full sub
  $2$-category of $\FFib/\cC$ consisting stacks over $\cC$  fibred in
  groupoids (resp. sets, resp. setoids).
 Note that in fact $\SStkS/\cC$ is a category.
 Then \eqref{eq_emb} can be extended into the following commutative
  diagram of $2$-categories
 \eq{ \label{eq_emb_ext}
 \xymatrix{
&\SSh(\cC) \ar[r]^-{\eta_\cC|_{\SSh(\cC)}}_-{\approx} \ar@{^(->}[d]
  &\SStkS/\cC\ar@{^(->}[r] \ar@{^(->}[d]
  &\SStkP/\cC\ar@{^(.>}[r] \ar@{^(.>}[d]
  &\SStk/\cC \ar@{^(.>}[d] \\
\cC \ar[r]^-{\epsilon_\cC}
   &\PPSh(\cC)\ar[r]^-{\eta_\cC}_-{\approx}
  &\FFibS/\cC\ \ar@{^(->}[r] &\FFibP/\cC\ \ar@{^(.>}[r] &\FFib/\cC
 }}
  where (dash) hook arrows indicate the $2$-functors defined by
  (strictly) full sub $2$-categories,
  $\SStkP/\cC$ (resp. $\FFibP/\cC$) is the essential images of the inclusion
  $2$-functor $\SStkS/\cC\hra \SStk/\cC$ (resp. $\FFibS/\cC\hra \FFib/\cC$).
 See \citesppp{042X}{0430}{0431}.
\item \label{it_comp_stk}
 In particular if $\cC=\cT_\tau$, where $\cT\ra\cS\in\FFib/\cS$, then
  the image of $\epsilon_{\cT}$ in \eqref{eq_emb_ext} is in $\SSh(\cT_\tau)$
  and  $\FFibR/\cT\hra \SStkP/\cT_\tau$.
 We also have $(\cY\ra \cT\ra \cS)\in \FFib/\cS$,
  whose topology inherited from
  $\cT_\tau$ agrees with $\cY_\tau$ (inherited directly from $\cS_\tau$).
 If moreover $\cT\ra\cS\in\SStk/\cS_\tau$,
  then $(\cY\ra \cT\ra \cS)\in \SStk/\cS_\tau$ (c.f. \citesp{09WX}).
\item A stack $\cY\ra \cC$ is a \emph{gerbe over $\cC$} if
 \enmt{[\upshape (i)]
 \item for any $U\in \cC$, there exists a covering $\{U_i\ra U\}$ in $\cC$
   such that $\cY_{U_i}$ is nonempty, and
 \item for any $U\in \cC$ and $x,y\in \cY_U$, there exists a covering
   $\{U_i\ra U\}$ in $\cC$ such that $x|_{U_i}\cong y|_{U_i}$ in $\cY_{U_i}$.
 }
 If $\cY\ra \cC$ has a section, we call $\cY$ \emph{trivial}.
\item \label{it_defi_Gerb}
 If a gerbe $\cY\ra \cC$ is \emph{bounded} by  $\sG\in\AAb(\cC)$
 \cite[IV.2.2]{giraud71cohnonab}, we denote it by $\cY\os{\sG}{\ra}\cC$.
 One may define \emph{$\sG$-equivalences} between them
  \cite[IV.2.2.7]{giraud71cohnonab}.
 In particular, it is an equivalence in $\SStk/\cS$.
 Denote by $\Gerb(\cC, \sG)$ the full sub $2$-category of $\SStk/\cC$ whose
  objects are gerbes over $\cC$ bounded by $\sG$.
 The set of $\sG$-equivalence classes in $\Gerb(\cC, \sG)$
  is denoted by $\sfH^2(\cC, \sG)$.
 Unless pointed out explicitly, we always consider gerbes bounded by
  \emph{abelian} sheaves $\sG\in\AAb(\cC)$.
 Then all trivial gerbes form a $\sG$-equivalence class,
  called the \emph{neutral} element, which
  makes $\sfH^2(\cC, \sG)$  a  pointed set.

}}

As an analogue to  descent by torsors (see \ref{prop_dbt}),
 we give
\thm{[Descent by gerbes] \label{thm_dbg}
Suppose that $\cA\in\SStkS/\cS_\tau$ and $\cX\in\SStk/\cS_\tau$.
Let $\sG\in\AAb(\cS_\tau)$ and  $f: \cY\ra \cX\in\Gerb(\cX_\tau, \sG)$.
Then $\Ob(\cXA^f)$ is characterized by
\eqn{
\cXA^f=\bigcup_{\s\in H_\tau^2(\cS, \sG)} f^\s(\cY^\s(\cA)),
}
 where $f^\s: \cY^\s\ra \cX\in\Gerb(\cX_\tau, \sG)$ is the twist of $\cY$ by
 $\s$ {\upshape (}see {\upshape \ref{defi_twist})}.
}
The proof will be given in \ref{pf_thm_dbg}.
Note that the assumption that
 $\cA\in\SStkS/\cS_\tau$ and $\cX\in\SStk/\cS_\tau$
 holds automatically in the
 classical  case where $X$ is a $k$-variety.
See \eqref{eq_emb_ext} and \ref{defi_gerbe} \eqref{it_comp_stk}.

\cor{ \label{cor_sdesc_by_dbg}
Suppose that $\cA\in\SStkS/\cS_\et$ and $\cA\in\SStk/\cS_\et$.
Then we have
\eqn{
\cXA^\sdesc=
 \bigcap_{\os{\sG\in\AAb(\cS_\et)}{f: \cY\ra \cX\in\Gerb(\cX_\et, \sG)}}
 \bigcup_{\s\in H_\et^2(\cS, \sG)} f^\s(\cY^\s(\cA)).
}}
\pf{
This follows immediately from \ref{ob_scoh} and \ref{thm_dbg}.
}

\cor{ \label{cor_rational_by_dbg}
With assumptions and notations in {\upshape \ref{thm_dbg}}, we have
\eqn{
\cXS=\bigcup_{\s\in H_\tau^2(\cS, \sG)} f^\s(\cY^\s(\cS)).
}}
\pf{
Take $\cA=\cS$ in \ref{thm_dbg}, and note that $\cXS^f=\cXS$.
Then the result follows.
}

\subsection{Classify gerbes} \label{gerb_H2}
We recall notations in \ref{coh}.
Let  $\tau\in\{\et, \fppf\}$ and  $\cT\in\FFib/\cS$.
Recall that $\cT_\tau$ is with the inherited topology from
 $\cS_\tau=(\SSch/S)_\tau$.
From now on, we shall write $\sfH_\tau^2(\cT, -)$ for $\sfH^2(\cT_\tau, -)$.
For  $\sG\in\AAb(\cT_\tau)$, there is a one-to-one correspondence
 (c.f.  \cite[IV.3.4.2 (i)]{giraud71cohnonab})
\eq{ \label{eq_gerb_H2}
\sfH_\tau^2(\cT, \sG) \os{\sim}{\ra} H_\tau^2(\cT, \sG),
}
 which maps the neutral element to $0$ and is functorial in $\sG$.

\subsection{The pull-back of a gerbe} \label{pb_gerb}
Let  $\tau\in\{\et, \fppf\}$ and $\sG\in\AAb(\cS_\tau)$ (whose pull-back to any
 $\cT_\tau$ is also denoted by $\sG$).
Throughout the rest of this section, we shall work in $\SStk/\cS_\tau$.
For any $1$-morphism $u: {\cT'}\ra\cT$ in $\SStk/\cS_\tau$ and
 $f: \cY\ra \cT\in\Gerb(\cT_\tau, \sG)$, we have the
\defi{[Pull-back of a gerbe] \label{defi_pb_gerb}
With the above notations,
 the \emph{pull-back} of the gerbe $\cY$ under $u$, denoted by $u^*\cY$, is the
 pull-back $u^*\cY$ as fibred categories (see \ref{defi_pb_fib}),
 indicated in the following $2$-cartesian diagram in $\FFib/\cS$
\eqn{\xymatrix{
u^*\cY=\cY\tm_\cT \cT'\ar[r]^-{f'}\ar[d]^-{u'} &\cT'\ar[d]^-u \\
\cY\ar[r]^-f &\cT
}}
Up to an equivalence, $f': u^*\cY\ra \cT'$ is a gerbe over $\cT'_\tau$ bounded
 by $\sG$ \citespp{06P2}{06P3},
 and the class $[u^*\cY]$ depends only on the class $[\cY]$.
}

Thus we have a well-defined map
\al{
u^*: \sfH_\tau^2(\cT, \sG) &\ra \sfH_\tau^2(\cT', \sG)  \label{eq_pb_gerb} \\
 [\cY] &\mpt [u^*\cY]. \nonumber
}
Under the correspondence \eqref{eq_gerb_H2} we have the following commutative
 diagram
\eq{ \label{eq_pb_corr}
\xymatrix{
\sfH_\tau^2(\cT, \sG)\ar[r]^-{u^*}\ar[d]^-\wr & \sfH_\tau^2(\cT', \sG)
 \ar[d]^-\wr \\
H_\tau^2(\cT, \sG)\ar[r]^-{u^*} &H_\tau^2(\cT', \sG)
}}
 where the lower map is defined by  \eqref{eq_pb_coh}.
See \cite[V.1.5.3]{giraud71cohnonab}.

\subsection{The contracted product of gerbes}
For $\sG\in\AAb(\cS_\tau)$ and $\cY,\cZ\in\Gerb(\cT_\tau, \sG)$,
 the \emph{contracted product} \cite[IV.2.4]{giraud71cohnonab} of $\cY$ and
 $\cZ$ exists, is also a gerbe over $\cT_\tau$ bounded by $\sG$,
 and is unique  up to a $\sG$-equivalence.
Denote it by $\cY\os{\sG}{\tm}\cZ$.
This gives a well-defined pairing
\aln{
(-, -)_c: \sfH_\tau^2(\cT, \sG)\tm  \sfH_\tau^2(\cT, \sG) &\ra
 \sfH_\tau^2(\cT, \sG) \\
 ([\cY], [\cZ]) &\mpt [\cY\os{\sG}{\tm}\cZ]
}
 which fits in to the following commutative diagram
\eq{ \label{eq_contr_corr}
\xymatrix{
\sfH_\tau^2(\cT, \sG)\tm \sfH_\tau^2(\cT, \sG) \ar[r]\ar[d]^-\wr
 &\sfH^2(\cT, \sG) \ar[d]^-\wr \\
H_\tau^2(\cT, \sG)\tm H_\tau^2(\cT, \sG) \ar[r] &H_\tau^2(\cT, \sG)
}}
 where the vertical maps come from the correspondence \eqref{eq_gerb_H2} and
 the lower pairing is the addition law of the
 abelian group $H_\tau^2(\cT, \sG)$.
See \cite[IV.3.3.2 (i), IV.3.4.2 (i)]{giraud71cohnonab}.

\lemm{ \label{lemm_pb_contr}
Let $u: \cT'\ra\cT$ be a $1$-morphism in $\SStk/\cS_\tau$.
Then we have a commutative diagram
\eqn{\xymatrix{
\sfH_\tau^2(\cT, \sG)\tm \sfH_\tau^2(\cT, \sG) \ar[r]^-{(-, -)_c}
 \ar[d]^-{u^*} &\sfH_\tau^2(\cT, \sG) \ar[d]^-{u^*} \\
\sfH_\tau^2(\cT', \sG)\tm \sfH_\tau^2(\cT', \sG)  \ar[r]^-{(-, -)_c}
 \ar[r] &\sfH_\tau^2(\cT', \sG)
}}}
\pf{
The diagram clearly comes from \eqref{eq_pb_corr}, \eqref{eq_contr_corr} and the
 fact that $u^*: H_\tau^2(\cT, \sG)\ra H_\tau^2(\cT', \sG)$
 is a homomorphism.
}

\subsection{} \label{twist}
Let
\eqn{\xymatrix@-2mm{
\cT'\ar[rr]^g \ar[dr]_-{a'} & &\cT\ar[dl]^-a \\
 &\cS
}}
 be a $1$-morphism in $\SStk/\cS_\tau$ and $\cZ\in\Gerb(\cS_\tau, \sG)$.
Then by the compositing property of $2$-cartesian diagrams,
 we have a natural $\sG$-equivalence
\eq{ \label{eq_pb_comp}
g^*\circ a^*\cZ\os{\approx}{\ra} a'^*\cZ.
}
\defi{[Twist of a gerbe] \label{defi_twist}
With notations in \ref{twist},
 let  $\s\in H_\tau^2(\cS, \sG)\os{\sim}{\ra} \sfH_\tau^2(\cS, \sG)$
 (see \eqref{eq_gerb_H2}), and $\cZ\ra \cS\in\Gerb(\cS_\tau, \sG)$
 represents  $\s$.
For any  $f: \cY\ra \cT\in  \Gerb(\cT_\tau, \sG)$,  we know the class
 $[\cY\os{\sG}{\tm}a^*\cZ]$ is independent of the choice of $\cZ$.
A \emph{twist} of $\cY$ by $\s\in H_\tau^2(\cS, \sG)$ is any  representative of
 the class $[\cY\os{\sG}{\tm}a^*\cZ]$,
 denoted by $f^\s: \cY^\s\ra\cT\in\Gerb(\cT_\tau, \sG)$.
Up to $\sG$-equivalences (which are equivalences in $\FFib/\cT$)
 $f^\s: \cY^\s\ra\cT\in\Gerb(\cT_\tau, \sG)$ is determined only by
 $\cY$ and $\s$.
}
\cor{ \label{cor_pb_twist}
With the above notations, we have
\eqn{
[(g^*\cY)^\s]=[g^*(\cY^\s)].
}}
\pf{
Considering \eqref{eq_pb_comp} and \ref{defi_twist}, one knows that
 this is a special case of  \ref{lemm_pb_contr}.
}

\subsection{Proof of $\ref{thm_dbg}.$ Theorem} \label{pf_thm_dbg}
By definition, $x\in \cXA^f$ means that the evaluation of
 $[\cY]\in H_\tau^2(\cX, \sG)$ at $x$ comes from $H_\tau^2(\cS,\sG)$, that is,
 there exists $\s\in H_\tau^2(\cS, \sG)$ such that
\eq{ \label{eq_Y(x)=-s}
[\cY](x)=-q^*\s \in H_\tau^2(\cA, \sG).
}
By \ref{cor_pb_twist}, we have $[(\cY_x)^\s]=[(\cY^\s)_x]$, where
 $\cY_x=x^*\cY$ and $(\cY^\s)_x=x^*(\cY^\s)$ are pull-backs
 (see \ref{defi_pb_gerb},  and note that $\cX,\cA\in\SStk/\cS_\tau$).
Since the   lower pairing in \eqref{eq_contr_corr} is in line with the addition
 law of the cohomology group, we also have
\eqn{
[\cY](x)+q^*\s=[\cY_x]+q^*\s=[(\cY_x)^\s].
}
Then \eqref{eq_Y(x)=-s} is equivalent to that
 $[(\cY^\s)_x]=0\in H_\tau^2(\cA, \sG)$, i.e., $(\cY^\s)_x$ is a trivial gerbe
 over $\cA_\tau$, which by definition means that $(\cY^\s)_x\ra \cA$ has a
 section in $\FFib/\cA$, i.e.,
 has a quasi-section in $\FFib/\cS$ since $\cA\in\FFibS/\cS$ (see \ref{rk_emb}
 \eqref{it_qs_s}).
By  \ref{lemm_sect_lift}, the last condition is to say that
 $x\in f^\s(\cY^\s(\cA))$. The proof is complete. \qed

\subsection{Non-abelian second descent obstruction} \label{sndesc}
With notations as in \ref{defi_gerbe} \eqref{it_defi_Gerb},
 we may also consider gerbes bounded by  (non-abelian) $\sG\in\SSGrp(\cC)$,
 and $\sG$-equivalences between them.
See \cite[IV.2.2, IV.2.2.7]{giraud71cohnonab}.
We keep the same notations such as  $\cY\os{\sG}{\ra}\cC$,  $\Gerb(\cC, \sG)$
 and  $\sfH^2(\cC, \sG)$.
\rk{
Note that
 if $\sG$ is not commutative, there may exists more than one neutral class in
 $\Gerb(\cC, \sG)$ and $\sfH^2(\cC, \sG)$ is not pointed in general.
}

Then   for $\sG\in\SSGrp(\cS_\tau)$,
 pull-backs of gerbes \ref{pb_gerb} also make sense.
Also \eqref{eq_pb_gerb} makes
 $\sfH_\tau^2(-, \sG): (\SStk/\cS_\tau)^\circ\ra \SSet$ a stable functor.
Thus we have  $\cXA^{\sfH_\tau^2(\cX, \sG)}$,
 the $\sfH_\tau^2(-, \sG)$-category, and the following
\defi{ \label{defi_sndesc}
Suppose that  $\cX,\cA\in\SStk/\cS_\et$.
Define the \emph{non-abelian second descent category} to be the full
 subcategory $\cXA^\sndesc$ of $\cXA$ whose objects are characterized by
\eqn{
\cXA^\sndesc=\bigcap_{\sG\in\SSGrp(\cS_\et)} \cXA^{\sfH_\et^2(\cX, \sG)}.
}}
From this, $\cXA^\sndesc$ is an obstruction category satisfying
 $\cXA^\sndesc\subseteq \cXA^\sdesc$. 

\section{The composite obstruction using second descent} \label{compob}
We now use results in previous two sections  to  give some full subcategories
 between $\cXS$ and $\cXA^\desc$.

\subsection{Composite obstruction categories}
Given an existing obstruction, we want to compose it with $\desc$ or $\sdesc$
 to produce new obstructions. To begin with, let us make the following
\defi{ \label{defi_func_map}
Let $\II\hra\FFib/\cS$ be a full sub $2$-category,
 $\ob$ a map sending each $\cX\in \II$ to an obstruction
 category (see \ref{defi_obcat}) of $\cXA$, called an \emph{obstruction map}
 on $\II$.
Define the \emph{trivial} obstruction map to be  the  constant map sending each
 $\cX\in \FFib/\cS$ to $\cXA$.

We say that an obstruction map  $\ob$ on $\II$
 is \emph{functorial} on $\II$ if for
 any $1$-morphism $f: \cY\ra \cX$ in $\II$, we have
 $f(\cYA^\ob)\subseteq \cXA^\ob$ as a full subcategory.

Note that if $\ob$ is functorial on $\II$, then so it is on any
 full sub $2$-category $\II'\hra\II$.
}
\lemm{ \label{lemm_Fob_func}
With the above notations,
 if $\ob=F: \II^\circ\ra \SSet$ is a stable functor
 $($see {\upshape \ref{defi_Fob}}$)$, then it is functorial on $\II$.

More generally, for a family of such stable functors $\{F_i\}_{i\in I}$,
 $\cXA^I:=\bigcap_{i\in I}\cXA^{F_i}$ is functorial on $\II$.
}
\pf{
For any $1$-morphism $f: \cY\ra \cX$ in $\II$, and $y\in\cYA^F$,
 since $\im F(y)\subseteq \im F(q)$, we have
\eqn{
\im F(f\circ y) = \im(F(y)\circ F(f)) \subseteq \im F(y)\subseteq \im F(q).
}
Then $f\circ y\in\cXA^F$ and hence $F$ is functorial on $\II$.

The second statement follows from the first obviously.
}
\coru{ \label{cor_*desc_func}
The maps $\desc$ (resp. $\sdesc$, resp $\sndesc$) is functorial on
 $\FFib/\cS$ (resp. $\FFib/\cS$, resp. $\SStk/\cS_\et$).
}

\prop{ \label{prop_func}
Let $\ob$ be a functorial obstruction map on $\II\hra\FFib/\cS$.
\enmt{[\upshape (i)]
\item \label{it_func_equiv}
 Let $\xymatrix{\cX'\ar@<1mm>[r]^f &\cX\ar@<1mm>[l]^g}$ be an equivalence in
  $\II$.
 Then we have $f(\cX'(\cA)^\ob)=\cXA^\ob$ and $\cX'(\cA)^\ob=g(\cXA^\ob)$.
\item \label{it_func_inv}
Let
\eqn{\xymatrix@-2mm{
\cY'\ar[rr]^-\theta \ar[dr]_-{f'} & &\cY\ar[dl]^-f \\
 &\cX
}}
 be a $2$-commutative diagram of $1$-morphisms in $\II$.
Then  $f'(\cY'(\cA)^\ob)\subseteq f(\cYA^\ob)$.
In particular, if $\theta$ is also an equivalence in $\FFib/\cX$,
 $($for example, $\cY'$ and $\cY$ are gerbes such that
 $[\cY']=[\cY]\in H_\tau^2(\cX, \sG)$$)$,
 then $f'(\cY'(\cA)^\ob)=f(\cYA^\ob)$.
}}
\pf{
Since $\ob$ is functorial, for \eqref{it_func_equiv},
 we have $\cX'(\cA)^\ob=g(f(\cX'(\cA)^\ob))\subseteq g(\cXA^\ob)\subseteq
 \cX'(\cA)^\ob$.
Thus $\cX'(\cA)^\ob=g(\cXA^\ob)$ and the other equality is similar.

For \eqref{it_func_inv}, the same reason yields
 $\theta(\cY'(\cA)^\ob)\subseteq \cYA^\ob$.
Then we have
\eqn{
f'(\cY'(\cA)^\ob)=f(\theta(\cY'(\cA)^\ob))\subseteq f(\cYA^\ob).
}
The proof is complete.
}

\defi{ \label{defi_deltaob}
For a functorial obstruction map $\ob$ on $\FFib/\cS$ (resp. $\SStk/\cS_\et$)
 we define the full subcategory $\cXA^\descob$ (resp. $\cXA^\sdescob$)
 whose objects are characterized by
\ga{
\cXA^\descob=
 \bigcap_{\os{\sG\in\SSGrp(\cS_\fppf)}{f: \cY\ra \cX\in\Tors(\cX_\fppf, \sG)}}
 \bigcup_{\s\in \check H_\fppf^1(\cS, \sG)} f^\s(\cY^\s(\cA)^\ob),
 \label{eq_descob}\\
(\text{resp. } \cXA^\sdescob=
 \bigcap_{\os{\sG\in\AAb(\cS_\et)}{f: \cY\ra \cX\in\Gerb(\cX_\et, \sG)}}
 \bigcup_{\s\in H_\et^2(\cS, \sG)} f^\s(\cY^\s(\cA)^\ob)). \label{eq_sdescob}
}}
\rk{ \label{rk_indep_twist}
These are well-defined. More precisely, by \ref{prop_func} \eqref{it_func_inv},
 $f^\s(\cY^\s(\cA)^\ob)$ is independent
 of the choice of the twisted gerbe (resp. torsor)
 $f^\s:\cY^\s\os{\sG}{\ra}\cX_\tau$ in
 its equivalence (resp. isomorphism) class.
}

\thm{ \label{thm_deltaob}
Let  $\delta\in\{\desc,\sdesc\}$.
Suppose that $\cA\in\FFibS/\cS$ $($resp. $\SStkS/\cS_\et)$ and
 $\II=\FFibP/\cS$ $($resp. $\SStk/\cS_\et)$ if
 $\delta=\desc$ $($resp. $\sdesc)$.
Let $\cX\in\II$ and $\ob$ be functorial on $\II$.
\enmt{[\upshape (i)]
\item \label{it_func_func} The category $\cXA^\deltaob$ is an obstruction
 category and $(\deltaob)$ is also functorial on $\II$.
\item \label{it_deltaob_in_cap} The objects of the category $\cXA^\deltaob$
 satisfies
\eqn{
\cXA^\deltaob\subseteq \cXA^\delta\cap \cXA^\ob.
}}}
\pf{
First note that in both cases, $\cX\in\II$ implies that  any
 $\cY\in\Tors(\cX_\fppf, \sG)$ (resp. $\Gerb(\cX_\et, \sG)$) is also
 an object in $\II$.
See \ref{cor_comp_fibp} (resp. \ref{defi_gerbe} \eqref{it_comp_stk}).
Moreover, $\II$ is closed under $2$-fibred products.
See \citespp{04SB}{02ZL}.

Since each  $\cY^\s(\cA)^\ob$ is an obstruction category,
 by  \ref{cor_sdesc_by_dbg} and \ref{cor_rational_by_dbg}
 (resp. \ref{cor_desc_by_dbt} and \ref{cor_rational_by_dbt}), we have
\eqn{
\cXS\subseteq \cXA^\deltaob\subseteq \cXA^\delta\subseteq \cXA.
}
In particular, $\cXA^\deltaob$ is an obstruction category.
Since $\ob$ is functorial on $\II$,
 $f^\s(\cY^\s(\cA)^\ob)\subseteq \cXA^\ob$.
Thus \eqref{it_deltaob_in_cap} is correct since $\cXA^\deltaob$ is an
 intersection of  unions of  $f^\s(\cY^\s(\cA)^\ob)$.

Next, we show that $(\deltaob)$ is functorial on $\II$.
For any $1$-morphism $g: \cX'\ra \cX$ in $\II$, we want to verify that
\eq{ \label{eq_deltaob_func}
g(\cX'(\cA)^\deltaob)\subseteq \cXA^\deltaob.
}
\enmt{
\item \label{it_delta=sdesc} $\delta=\sdesc$.
 Let $f: \cY\ra\cX_\et\in\Gerb(\cX_\et, \sG)$.
 We have the pull-back of $\cY$ under $g$ (see  \ref{defi_pb_gerb})
  $f': \cY':=g^*\cY\ra \cX'\in\Gerb(\cX'_\et, \sG)$.
 Let $\s\in H_\et^2(\cS, \sG)$ and
 \gan{
 f^\s: \cY^\s\ra\cX\in\Gerb(\cX_\et, \sG), \\
 f'^\s: \cY'^\s\ra\cX'\in\Gerb(\cX'_\et, \sG),
}
  be the corresponding twists (see \ref{defi_twist}).
 The pull-back $\cY''_\s:=g^*(\cY^\s)$ fits into the following
  $2$-cartesian diagram in $\SStk/\cS_\et$
 \eqn{\xymatrix{
 \cY''_\s\ar[r]^-{f''_\s}\ar[d]^-{g''_\s} &\cX'\ar[d]^-g \\
 \cY^\s\ar[r]^-{f^\s} &\cX
 }}
 Since $\ob$ is functorial on $\SStk/\cS_\et$, we have
 \eqn{
 g(f''_\s(\cY''_\s(\cA)^\ob))=f^\s(g''_\s(\cY''_\s(\cA)^\ob))\subseteq
  f^\s(\cY^\s(\cA)^\ob).
 }
 On the other hand, by  \ref{cor_pb_twist} we have
 \eqn{
 [\cY''_\s]=[g^*(\cY^\s)]=[(g^*\cY)^\s]=[\cY'^\s].
 }
 Thus by  \ref{prop_func} \eqref{it_func_inv},
  $f'^\s(\cY'^\s(\cA)^\ob)=f''_\s(\cY''_\s(\cA)^\ob)$.
 It follows that
 \aln{
 g(\bigcup_{\s\in H_\et^2(\cS,\sG)} f'^\s(\cY'^\s(\cA)^\ob)) &=
  g(\bigcup_{\s\in H_\et^2(\cS,\sG)} f''_\s(\cY''_\s(\cA)^\ob)) \\
  &\subseteq \bigcup_{\s\in H_\et^2(\cS,\sG)} f^\s(\cY^\s(\cA)^\ob).
 }
 If $f:\cY\ra \cX$ runs over $\Gerb(\cX_\et, \sG)$, then the resulting
  $f':\cY'\ra \cX'$ forms a subclass of $\Gerb(\cX'_\et, \sG)$.
 This proves \eqref{eq_deltaob_func}, i.e.,  $(\sdescob)$ is functorial
  on $\SStk/\cS_\et$.
\item $\delta=\desc$.
 The argument for this case is the same as \eqref{it_delta=sdesc}, except a few
  points.
 Let $\sY\in\Tors(\cX_\fppf, \sG)$ and $f: \cY\ra\cX\in\FFib/\cX$ be associated
  to it (see \ref{tors_as_fib}).
 We have the pull-back $g^*\sY\in\Tors(\cX'_\fppf, \sG)$ (see \ref{pb_tors}) and
  its associated $f': \cY'\ra \cX'\in\FFib/\cX'$.
 Let $\s\in \check H_\fppf^1(\cS, \sG)$ and
 \gan{
 f^\s: \cY^\s\ra\cX\in\Tors(\cX_\fppf, \sG^\s), \\
 f'^\s: \cY'^\s\ra\cX'\in\Tors(\cX'_\fppf, \sG^\s),
 }
  be the corresponding twists (see \ref{defi_twist_tors}).
 By \ref{lemm_pb_twist_tors} we have an isomorphism in $\FFib/\cX'$
 \eqn{\xymatrix@-2mm{
 \cY'^\s\ar[rr]^-\beta_-{\sim} \ar[dr]_-{f'^\s}
  & &\tilde\cY_\s\ar[dl]^-{\tilde h_\s} \\
  &\cX'
 }}
  where
  $(\tilde h_\s: \tilde\cY_\s\ra \cX'):=\eta_{\cX'}(g^*(\sY^\s))\in\FFib/\cS$.
 On the other hand, the pull-back $\cY''_\s:=g^*(\cY^\s)$
  (see \ref{defi_pb_fib}) fits into the
  following $2$-cartesian diagram in $\FFibP/\cS$
 \eqn{\xymatrix{
 \cY''_\s\ar[r]^-{f''_\s}\ar[d]^-{g''_\s} &\cX'\ar[d]^-g \\
 \cY^\s\ar[r]^-{f^\s} &\cX
 }}
 Since $\ob$ is functorial on $\FFibP/\cS$, we have
 \eqn{
 g(f''_\s(\cY''_\s(\cA)^\ob))=f^\s(g''_\s(\cY''_\s(\cA)^\ob))\subseteq
  f^\s(\cY^\s(\cA)^\ob).
 }
 We know from the proof of \ref{lemm_pb_eta} that there exists a $2$-commutative
  diagram in $\FFib/\cS$
 \eqn{\xymatrix@-2mm{
 \tilde\cY_\s\ar[rr]^-\alpha \ar[dr]_-{\tilde h_\s}
  & &\cY''_\s\ar[dl]^-{f''_\s} \\
  &\cX'
 }}
 Thus apply  \ref{prop_func} \eqref{it_func_inv} to $\beta$ and $\alpha$,
  we have
 \eqn{
 f'^\s(\cY'^\s(\cA)^\ob)=\tilde h_\s(\tilde\cY_\s(\cA)^\ob)\subseteq
  f''_\s(\cY''_\s(\cA)^\ob).
 }
 It follows that
 \aln{
 g(\bigcup_{\s\in \check H_\fppf^1(\cS,\sG)} f'^\s(\cY'^\s(\cA)^\ob))
  &\subseteq
  g(\bigcup_{\s\in \check H_\fppf^1(\cS,\sG)} f''_\s(\cY''_\s(\cA)^\ob)) \\
  &\subseteq \bigcup_{\s\in \check H_\fppf^1(\cS,\sG)} f^\s(\cY^\s(\cA)^\ob).
 }
 If $f:\cY\ra \cX$ runs over $\Tors(\cX_\fppf, \sG)$, then the resulting
  $f':\cY'\ra \cX'$ forms a subclass of $\Tors(\cX'_\fppf, \sG)$.
 This proves \eqref{eq_deltaob_func}, i.e.,  $(\descob)$ is functorial
  on $\FFibP/\cS$.
}
The proof is complete.
}
\rk{
If $\cA\in\SStkS/\cS_\et$ and $\cX\in\SStkP/\cS_\et$ (for example,
 $\cA, \cX\in\cS$ are schemes), then the assumptions
 in \ref{thm_deltaob} always hold.
}

By recursively using \ref{defi_deltaob} and
 \ref{thm_deltaob} \eqref{it_func_func} one can make the following
\defi{ \label{defi_comp_ob}
With notations and assumptions in \ref{thm_deltaob},
 for $n\ge1$, recursively define the  obstruction map
 $(\delta^n,\ob)=(\delta, (\delta^{n-1}, \ob))$
 and $(\delta^1,\ob)=(\delta,\ob)$.
}
\cor{ \label{cor_comp_ob}
With notations and assumptions in {\upshape \ref{defi_comp_ob}}, we have
\eqn{
\cXS\subseteq \dots\subseteq \cXA^{\delta^{n+1},\ob}\subseteq
 \cXA^{\delta^n,\ob}\subseteq \dots \subseteq\cXA^{\delta^1,\ob}\subseteq
 \cXA^\delta\cap \cXA^\ob\subseteq \cXA
}
 and that all obstruction maps above are functorial on $\II$.
In particular, we may take $\ob=\desc$, $\sdesc$ or the trivial obstruction.
We may also take $\ob=\sndesc$ but restrict $\delta=\sdesc$ only.
}
\pf{
In  \ref{thm_deltaob}, recursively use \eqref{it_func_func}, and
 the inclusions (as full subcategories) follows from \eqref{it_deltaob_in_cap}.
}
\rk{
We have some remarks on  \ref{cor_comp_ob}.
\enmt{[\upshape (i)]
\item In next section,
  we shall see that we my also take $\ob=\hdesc{i}{\tau}$ and  $\FF$
  in \ref{cor_comp_ob}.
\item From \ref{cor_comp_ob} we know that
  $\cXA^{\delta^n,\desc}\subseteq \cXA^\desc$ and
  $\cXA^{\desc^i,\ob}\subseteq \cXA^\desc$.
 Combined with \ref{prop_cXAH1}, we note that in the classical case where
  $X$ is a $k$-variety, obstruction maps  $(\delta^n, \desc)$, $(\desc^i, \ob)$
  yield obstruction sets contained in the classical descent set $\XA^\cdesc$.
}}

\section{Higher descent obstruction and derived obstruction} \label{hdescder}
We propose some new kinds of obstructions in this section and exploit some of
 their properties.

\subsection{Higher descent categories} \label{hdesc}
In  \ref{Hiob} we obtain  $\cXA^{H_\tau^i(\cX, \sG)}$,
 the $H_\tau^i(-, \sG)$-category, where  $\sG\in\AAb(\cS_\tau)$ and $\tau\in
 \{\et,\fppf\}$.
\defi{ \label{defi_idesc}
For $i\ge0$, define the \emph{$i$-th $\tau$-descent category} to be the full
 subcategory  $\cXA^\hdesc{i}{\tau}$ of $\cXA$ whose
 objects are characterized by
\eqn{
\cXA^\hdesc{i}{\tau}=\bigcap_{\sG\in\AAb(\cS_\tau)} \cXA^{H_\tau^i(\cX, \sG)}.
}}
By \ref{lemm_Fob_func}, $\hdesc{i}{\tau}$ is functorial on $\FFib/\cS$.
Thus we may take $\ob=\hdesc{i}{\tau}$ in \ref{cor_comp_ob}.

It is clear that $\sdesc=\hdesc{2}{\et}$ and for $\desc$ we have the following
\prop{
If $\cX,\cA\in\FFibR/\cS$, then
 $\cXA^\desc\subseteq \cXA^\hdesc{1}{\fppf}$.
}
\pf{
This immediately follows from \ref{defi_idesc} and \ref{prop_cXAH1}
 \eqref{it_cXAH1_ab}.
}

\prop{ \label{prop_hdesc_cmp}
For $i\ge1$ and $\tau\in\{\et, \fppf\}$, we have
\eqn{
\cXA^\hdesc{i}{\tau}\subseteq \cXA^\hdesc{(i+1)}{\tau}.
}}
\pf{
Let $\sG\in\AAb(\cS_\tau)$, and $\sF^.\in C(\cS_\tau)^{\ge0}$ be
 an injective resolution of $\sG$ where $C(\cS_\tau)$ is the category of
 complexes of sheaves in $\AAb(\cS_\tau)$.
Then we have a short exact sequence
\eqn{
0\ra \sG\ra \sF\ra \sH\ra 0,
}
 where $\sF=\sF^0$ and $\sH=\tau_{\le1}\sF^.$.
Since $\sF$ is injective and the pull-back of sheaves
 is an exact functor and sends acyclic sheaves
 to acyclic ones,
 it follows that for each $i\ge1$, we have a natural isomorphism of functors
\eqn{
\delta^i: H_\tau^i(-, \sH)\os{\sim}{\Ra} H_\tau^{i+1}(-, \sG):
 (\FFib/\cS)^\circ\ra \AAb.
}
Then by \ref{lemm_ob_nat} below we have
 $\cX^{H_\tau^i(\cX, \sH)}=\cXA^{H_\tau^{i+1}(\cX, \sG)}$.
Thus let $\sG$ runs over all objects in $\AAb(\cS_\tau)$, and we obtain the
 desired inclusion.
}

\lemm{ \label{lemm_ob_nat}
Let $F, G: (\FFib/\cS)^\circ\ra \SSet$ be two  stable functors and
 $\xi: F\Ra G$  a natural transformation.
For $A\in F(\cX)$ and $B=\xi_\cX(A)$, we have  $\cXA^A\subseteq \cXA^B$.

In particular, if $\xi$ is a natural isomorphism, $\cXA^F=\cXA^G$.
}
\pf{
The result clearly follows from the commutative diagram
\eqn{\xymatrix@-2mm{
\cXS\ar[rr]\ar[dd]_-{B(-)}\ar[dr]^-{A(-)} &&\cXA\ar[dd]_(.3){B(-)}
 \ar[dr]^-{A(-)} \\
 &F(\cS)\ar@{-}'[r]\ar[dl]^-{\xi_\cS}
 & &F(\cA)\ar[dl]^-{\xi_\cA}\ar@{<-}'[l] \\
G(\cS)\ar[rr] &&G(\cA)
}}}

\subsection{} \label{prod}
Besides the obstructions constructed from functors
 (see \ref{Fob}),
 Harpaz and Schlank \cite{hs13homotopy} developed a new kind of obstruction
 (\emph{homotopy obstruction})
 using \'etale homotopy theory, and established
 its connections with some existing obstructions such as $\Br$, $(\etBr)$ and
 $\cdesc$.
One consequence is the following
\thmu{[{\cite[Thm. 9.147]{hs13homotopy}}]
Let $X$ and $Y$ be smooth geometrically connected varieties over a number field
 $k$. Then we have
\eqn{
(X\tm_k Y)(\bfA_k)^\etBr=\XA^\etBr\tm \YA^\etBr.
}}
This result ensured the commutativity
 of  Brauer-Manin sets with taking product of two varieties.
The same result with $(\etBr)$
 replaced by $\Br$ was given by Skorobogatov and Zarhin \cite{sz14product}
 for smooth projective geometrically integral varieties  and
 the author \cite{lv20brprodge} for open varieties.
Next we  propose another new kind of obstructions for fibred categories,
 given by pseudofunctors.
In particular we define the derived obstruction which also has good behavior
 under a product and,
 with topology chosen properly,
 is not larger than $\desc$, $\sdesc$ and $\hdesc{i}{\tau}$.

\subsection{Obstructions given by pseudofunctors} \label{psfob}
Let $\FF: \FFib/\cS\ra \CC$ be a pseudofunctor to a $2$-category $\CC$.
Then we have the following $2$-commutative diagram in $\CCat$
\eqn{\xymatrix{
\cXS\ar[r]^-{\MMor_{\FFib/\cS}(q, \cX)} \drtwocell\omit{_\gamma}
 \ar[d]^-{\FF} &\cXA\ar[d]^-{\FF} \\
\MMor_\CC(\FF(\cS), \FF(\cX))\ar[r]^{\varphi_\FF(\cX)}
  &\MMor_\CC(\FF(\cA), \FF(\cX))
}}
 where $\varphi_\FF(\cX)=\MMor_\CC(\FF(q), \FF(\cX))$ and
 $\gamma$ is constructed by the definition of pseudofunctors.
Explicitly, for each $z\in\cXS$, $\gamma_z$ is the
 $2$-isomorphism in $\CC$
\eq{ \label{eq_n_nat}
\gamma_z: (\FF\circ\MMor_{\FFib/\cS}(q, \cX))(z)=
 \FF(z\circ q)\os{\sim}{\Ra} (\varphi_\FF(\cX)\circ\FF)(z)=\FF(z)\circ\FF(q).
}
\defi{
The $\FF$-\emph{category} is  the full
 subcategory $\cXA^\FF$ of $\cXA$ whose objects are characterized by
\eqn{
\cXA^\FF=\{x\in \cXA\mid \text{ $\FF(x)$ is in the essential image of
 $\varphi_\FF(\cX)$} \}.
}}
By \eqref{eq_n_nat}, one checks that $\cXA^\FF$ is also an obstruction category
 (see \ref{defi_obcat}), and by an argument analogue to \ref{lemm_Fob_func},
 we know that $\FF$ is functorial on $\FFib/\cS$.
Thus we may take $\ob=\FF$ in  \ref{cor_comp_ob}.

\subsection{Modified $2$-category of $\St{\cS}$-topoi} \label{mod_tps}
For any $\cT\in\FFib/\cS$,
 the structure morphism $\cT\ra \cS$ induces a morphism of topoi
 $\St{\cT}\ra \St{\cS}$  (c.f. \ref{pb_sh}), which makes $\St{\cT}$ a
 $\St{\cS}$-topos.
Consider the $2$-category of $\St{\cS}$-topoi,
 we use  \cite[VIII.0]{giraud71cohnonab} as  whose definition,
 but a slight modification as follows.
Objects are still morphism of topoi $x: X\ra\St{\cS}$ and
 we modify $\MMor(X, Y)$ into a full subcategory, so that,
 in the definition of a $1$-morphism
\eqn{\xymatrix@-2mm{
X\ar[rr]^-f \drlowertwocell<0>_{x}{<-2>_\alpha}
 & &Y\ar[dl]^-y \\
 &\St{\cS}
}}
 which is a $2$-isomorphism of topoi $\alpha: yf\os{\sim}{\Ra}x$,
 we further require that $\alpha^*=\id_{x^*}: x^*=f^*y^*$.
This is the case if the $1$-morphism comes from the one in $\FFib/\cS$.
Denote by    $\TTpsSt$
 the \emph{modified} $2$-category of $\St{\cS}$-topoi above.
\defi{[Derived obstruction] \label{defi_do}
With notations in \ref{psfob}, let $\CC=\TTpsSt$ and
 $\FF: \FFib/\cS\ra \TTpsSt$ be the pseudofunctor defined by
\eqn{
\cT\mpt \St{\cT}, \quad f\mpt f=(f^*,f_*),
 \quad \xi\mpt \xi=(\xi_*,\xi^*).
}
This is well-defined by  \ref{pb_sh} and the above definition of
 $\TTpsSt$ \ref{mod_tps}.
Then we define the \emph{derived $($obstruction$)$ category}
 to be the category defined by $\FF$ given above, denoted by $\cXA^\SSht$.
In other words, $x\in\cXA^\SSht$ if and only if there exists a $1$-morphism
 $x_0: \St{\cS}\ra \St{\cX}$ along with a $2$-isomorphism
 $x\os{\sim}{\Ra} x_0 q$ in $\TTpsSt$.
As all $\FF$ do,  $\SSht$ is functorial on $\FFib/\cS$.
}

\subsection{} \label{ops}
Let $\La$ be a commutative ring and $\cT\in\FFib/\cS$.
Denote  by $\MMod(\cT_\tau, \La)$ the sheaves of $\La$-modules on $\cT_\tau$,
 and by $D(\cT_\tau, \La)$ the corresponding derived category.
In particular, we have $\AAb(\cT_\tau)=\MMod(\cT_\tau, \ZZ)$ and
 $D(\cT_\tau)=D(\cT_\tau, \ZZ)$.
Now we suppose that  $x\in \cXA^\SSht$. Thus we have a $2$-isomorphism
 $x\os{\sim}{\Ra} x_0 q$ in $\TTpsSt$ for some $x_0$.
Then we have isomorphisms of functors
\ga{
x_*\os{\sim}{\Ra} x_{0*}q_*: \SSh(\cA_\tau)\ra \SSh(\cX_\tau),\
 \SSGrp(\cA_\tau)\ra \SSGrp(\cX_\tau) \text{ or }
 \MMod(\cA_\tau, \La)\ra \MMod(\cX_\tau, \La) \text{, and} \nonumber \\
x^*\os{\sim}{\Ra} q^*x_0^*: \SSh(\cX_\tau)\ra \SSh(\cA_\tau),\
 \SSGrp(\cX_\tau)\ra \SSGrp(\cA_\tau) \text{ or }
 \MMod(\cX_\tau, \La)\ra \MMod(\cA_\tau, \La), \label{eq_do_Sh}
}
 which also induce isomorphisms of derived functors \cite[18.6]{ks06cat}
\ga{
Rx_*\os{\sim}{\Ra} Rx_{0*}Rq_*: D(\cA_\tau, \La)\ra D(\cX_\tau, \La)
 \text{ and} \nonumber \\
x^*\os{\sim}{\Ra} q^*x_0^*: D(\cX_\tau, \La)\ra D(\cA_\tau, \La).
 \label{eq_do_D}
}

\prop{ \label{prop_do_cmp}
We have  $\cXA^{\SSh_\fppf}\subseteq \cXA^\desc$  and
 $\cXA^\SSht\subseteq \cXA^\hdesc{i}{\tau}$ for all $i\ge0$.
}
\pf{
For any $1$-morphism  $f: {\cT'}\ra\cT$ in $\FFib/\cS$ and $L\in D(\cT_\tau)$,
 recall the functorial commutative diagram \eqref{eq_pb_coh_tps}.
In particular, for $\sG\in\AAb(\cS_\tau)$, we have an isomorphism of functors
\eqn{
H_\tau^i(-, \sG)\os{\sim}{\Ra} \Hom_{D(-)}(\ZZ, \sG[i]).
}
Then the second inclusion follows from  \eqref{eq_do_D}, \ref{mod_tps} and
 the definition of $\cXA^\hdesc{i}{\tau}$ (\ref{defi_idesc}).

For the first inclusion, note that \eqref{eq_pb_tors_gen}, in particular
 \eqref{eq_pb_tors} still makes sense for any $1$-morphism of topoi
 $u: \St{\cT'}\ra \St{\cT}$.
Take $\tau=\fppf$ and
 the result follows from \eqref{eq_do_Sh}, \ref{mod_tps} and
 the definitions of $\cXA^\desc$ (\ref{defi_desc}).
}

Apply \ref{cor_comp_ob} and \ref{prop_do_cmp} we obtain the following
\coru{ \label{cor_comp_do}
With notations and assumptions in \ref{defi_comp_ob},
 write $\delta_\et=\sdesc$ and $\delta_\fppf=\desc$.
Then  for any $n\ge1$ we  have
\eqn{
\cXS\subseteq \cXA^{\delta^n,\SSht}\subseteq
 \cXA^{\delta^1,\SSht}\subseteq \cXA^\delta\cap\cXA^\SSht\subseteq
 \cXA^\SSht \subseteq \cXA^{\delta_\tau}\subseteq \cXA,
}
 and that all obstruction maps above are functorial on $\II$.
}
\rk{
Thus $\cXA^{\delta^n,\SSht}$ with $n$ and $\tau$ chosen properly, is
 not larger than
 all obstructions occurring in this text, except $\cXA^{\sdesc^n,\sndesc}$.
}

\subsection{} Now we consider the  behavior of derived obstructions
 under a product.
Recall that a $1$-morphism $f: \cY\ra \cX$ in $\FFib/\cS$  is called
 \emph{representable} if for every $1$-morphism $U\ra \cX$ with $U\in\cS$,
 the $2$-fibred product $U\tm_\cX \cY$ is representable (in $\FFib/U$).
\thm{ \label{thm_do_prod}
Consider the following $2$-cartesian diagram in $\FFib/\cS$
\eqn{
\xymatrix{
\cV\ar[r]^-{j'}\ar[d]^-{f'} &\cY\ar[d]^-f \\
\cU\ar[r]^-j &\cX
}}
 where $\cU$ and $f$ are representable
 $($hence $\cV$ is necessarily representable$)$.
Then we have a $2$-cartesian diagram in $\CCat$
\eqn{
\xymatrix{
\cVA^\SSht\ar[r]^-{j'}\ar[d]^-{f'} &\cYA^\SSht\ar[d]^-f \\
\cUA^\SSht\ar[r]^-j &\cXA^\SSht
}}
 that is,
 projections to $\cU$ and $\cX$ induce a natural equivalence in $\CCat$
\eqn{
\cVA^\SSht\os{\approx}{\lra} \cUA^\SSht\tm_{\cXA^\SSht} \cYA^\SSht.
}}
\pf{
By the definition of $2$-cartesian diagrams,
 projections to $\cU$ and $\cX$ induce a natural equivalence
\eqn{
\varpi: \cVA\os{\approx}{\lra} \cUA\tm_{\cXA} \cYA.
}
Since $\SSht$ is fuctorial, and the natural functor
\eqn{
\cUA^\SSht\tm_{\cXA^\SSht} \cYA^\SSht\lra \cUA\tm_{\cXA} \cYA
}
 identifies  $\cUA^\SSht\tm_{\cXA^\SSht} \cYA^\SSht$ as a full subcategory of
 $\cUA\tm_{\cXA} \cYA$,
 we see that
 $\varpi$ restricts to a fully faithful
 functor
\eqn{
\varpi^\SSht: \cVA^\SSht\lra \cUA^\SSht\tm_{\cXA^\SSht} \cYA^\SSht
}
 which is the functor in the statement.
Then it suffices to show that $\varpi^\SSht$ is essentially full.
Let $t=(u, y, \alpha)\in \cUA^\SSht\tm_{\cXA^\SSht} \cYA^\SSht$, where
 $u\in\cU^\SSht$, $y\in\cY^\SSht$ and $\alpha: ju\os{\sim}{\ra} fy$
 is an isomorphism in $\cXA^\SSht$.
Find $v\in\cVA$ such that $\varpi v\cong t$.
Thus we have the following  $2$-commutative diagram in $\FFib/\cS$
\eqn{\xymatrix{
\cA\ar@/^/[drr]^-{y} \ar[dr]|-{v} \ar@/_/[ddr]_-{u} \\
 &\cV\ar[r]^-{j'}\ar[d]_-{f'} &\cY\ar[d]^-f \\
 &\cU\ar[r]^-j &\cX
}}
 which gives the corresponding $2$-commutative diagram in $\TTps/\St{\cS}$
\eq{ \label{eq_tps}
\xymatrix{
\St{\cA}\ar@/^/[drr]^-{y} \ar[dr]|-{v} \ar@/_/[ddr]_-{u} \\
 &\St{\cV}\ar[r]^-{j'}\ar[d]^-{f'}
 &\St{\cY}\ar[d]^-f \\
 &\St{\cU}\ar[r]^-j &\St{\cX}
}}
Since $\cU$ and $j$ are representable, we reduce to the case that
 $\cU=U, \cV=V\in \cS$ are schemes.
Then by the $2$-Yoneda lemma, the morphism $j: U\ra \cX$ accounts to an object
 $u'\in \cX_U$.
It follows that we may assume $j$ is the embedding $\cX/u'\ra \cX$ and
 similarly $j': \cY/y'\ra \cY$ where $y'\in \cY_V$.

In fact, since $f$ is representable, there exists a continuous functor
 $h: \cX_\tau\ra \cY_\tau$ which induces the same morphism of topoi
 $\St{\cY}\ra \St{\cX}$ as $f$ does, and we have $y'=hu'$.
See \citesp{06W8}.
Thus we may assume the $2$-square in \eqref{eq_tps} is of the form
\eqn{
\xymatrix{
\St{\cY}/\ul{y'}\ar[r]^-{j_{\ul{y'}}}\ar[d]^-{f|_{\ul{y'}}} &\St{\cY}\ar[d]^-f \\
\St{\cX}/\ul{u'}\ar[r]^-{j_{\ul{u'}}} &\St{\cX}
}}
 which is $2$-cartesian in $\TTpsSt$ by \cite[VIII.1.6]{giraud71cohnonab}
 (one checks that  this $2$-diagram is $2$-cartesian not only  in the original
 $2$-category but also  in the modified one).
Thus the $2$-square in \ref{eq_tps} is also $2$-cartesian.

Since $u\in \cU^\SSht$ and $y\in \cY^\SSht$, there exist morphisms
 $u_0: \St{\cS}\ra \St{\cU}$ and
 $y_0: \St{\cS}\ra \St{\cY}$ such that
 we have $2$-isomorphisms
 $u\os{\sim}{\Ra} u_0 q$ and
 $y\os{\sim}{\Ra} y_0 q$.
Then  the $2$-cartesian square in \ref{eq_tps} ensures that
 there exists a $1$-morphism
 $v_0: \St{\cS}\ra \St{\cV}$ with $2$-isomorphisms
\eq{ \label{eq_v_0}
f'v_0\os{\sim}{\Ra} u_0 \text{ and }
 j'v_0\os{\sim}{\Ra} y_0.
}
Then by universal property, \eqref{eq_v_0},
 and the $2$-diagram \eqref{eq_tps} with the square $2$-cartesian,
 there exists a $2$-isomorphism $v\os{\sim}{\Ra} v_0q$.
It follows that $v\in \cV^\SSht$ and $\varpi^\SSht v\cong t$.
The proof is complete.
}
\cor{ \label{cor_do_prod_over_cS}
With notations and assumptions as in Theorem \ref{thm_do_prod}, if $\cX=\cS$,
 then projections to $\cU$ and $\cX$ induce a natural equivalence
\eqn{
\cVA^\SSht\os{\approx}{\lra} \cUA^\SSht\tm \cYA^\SSht.
}}
\pf{
Just note that if $\cX=\cS$, then
 $\cS(\cA)^\SSht=\cS(\cS)=\{q\}$ is a singleton category.
}
\coru{ \label{cor_do_prod_classical}
Let $X$ and $Y$ be $k$-schemes. Then we have
\eqn{
(X\tm_k Y)(\bfA_k)^\SSht=\XA^\SSht\tm \YA^\SSht.
}}
In particular, taking $\tau=\fppf$ and in view of \ref{prop_do_cmp} and
 \eqref{eq_desc_in_cdesc},
 we obtain  the same product formulate as mentioned  in \ref{prod} but with
 the derived obstruction $\cXA^{\SSh_\fppf}$ not larger than all known
 obstructions
 before this text.

\rk{ \label{rk_do_bit_small}
If we use small \'etale site for derived obstructions, things will get in
 trouble:
\enmt{[\upshape (i)]
\item For schemes and algebraic spaces, the pseudofunctor
  $X\mpt (X_\et, \cO_X)$ assigning  each $X$ to its locally ringed small
  \'etale topos is fully faithful (morphisms up to $2$-isomorphisms on the
  right  hand side). See \citespp{04I7}{04KL}.
 Thus derived obstructions using this pseudofunctor lead to rational points
  and make no sense.
\item For \DM stacks with the small \'etale site topology, in order to produce
  the same  kind of cartesian square of topoi, the
  morphism $j$ in the cartesian square in Theorem \ref{thm_do_prod} should be
  representable \'etale \cite[Constructions 2.7, 2.8]{zheng15six},
  which is too strong in the classical case.
 For example in  \ref{cor_do_prod_classical} we should require
  that $X$ is an \'etale variety over $\bfA_k$.

 In general, the $2$-squares of topoi come from $2$-cartesian squares in
  $\FFib/\cS$ are  not necessarily $2$-cartesian.
}}

\section*{Acknowledgment}
The author would like to thank the referees for valuable suggestions, and
 Junchao Shentu, Hang Yin for  helpful discussions.

\bibliography{\bibfilename}
\bibliographystyle{amsplain}
\end{document}